\input amstex
\input epsf
\documentstyle{amsppt}
\magnification=1200
\pretolerance=200
\tolerance=400
\TagsOnRight
\NoRunningHeads
\define\({\left(}
\define\){\right)}
\define\a{\alpha}
\define\A{\Cal A}
\define\amb{\allowmathbreak}
\define\dd{\partial}
\define\e{\varepsilon}
\define\emb{\hookrightarrow}
\define\f{\varphi}
\define\G{\Gamma}
\define\GL{\operatorname{GL}}
\define\id{\operatorname{id}}
\redefine\Im{\operatorname{Im}}
\define\Int{\operatorname{Int}}
\define\lk{\operatorname{lk}}
\define\oa{\overarrow}
\define\ol{\overline}
\define\ph{\phantom}
\define\pusto{\varnothing}
\define\R{\Bbb R}
\define\sign{\operatorname{sign}}
\define\sk{\operatorname{sk}}
\define\SL{\operatorname{SL}}
\define\sm{\setminus}
\define\T{\operatorname{T}}
\define\th{\theta}
\define\thc{$\theta$-curve\ }
\define\thcs{$\theta$-curves\ }
\define\tM{\widetilde M}
\define\tr{\operatorname{tr}}
\define\tto{\longrightarrow}
\define\wt{\widetilde}
\define\x{\times}
\define\Z{\Bbb Z}


\topmatter
\title Towards Lower Bounds for Complexity of 3-Manifolds: a Program \endtitle
\author S.~Anisov \endauthor
\address Institute for System Studies, Russian Academy of Sciences \endaddress
\curraddr IHES, Le Bois-Marie, 35, Route de Chartres, F-91440 Bures-sur-Yvette,
France \endcurraddr
\email anisov\@mccme.ru \endemail
\thanks This research was supported by FCT-Portugal through the Research Units
Pluriannual Funding Program and by EPDI (European Post-Doctoral Institute).
\endthanks
\abstract For a 3-dimensional manifold $M^3$, its complexity $c(M^3)$,
introduced by S.~Matveev, is the minimal number of vertices of an almost simple
spine of~$M^3$; in many cases it is equal to the minimal number of tetrahedra
in a singular triangulation of~$M^3$. An approach to estimating $c(M^3)$ from
below for total spaces of torus bundles over~$S^1$, based on the study of \thcs
in the fibers, is developed, and pseudominimal special spines for these
manifolds are constructed, which we conjecture to be their minimal spines. We
also show how to apply some of these ideas to other 3-manifolds.   \endabstract
\toc\nofrills{Contents}
\widestnumber\head{\S9.9.}
\specialhead{\S 1.} Introduction			\endspecialhead
    \head{1.1.} Brief outlook		\page2		\endhead
    \head{1.2.} Definitions		\page3		\endhead
    \head{1.3.} Example: lens spaces	\page5		\endhead
\specialhead{\S 2.} Upper bound				\endspecialhead
    \head{2.1.} $\th$-curves		\page7		\endhead
    \head{2.2.} On $\SL(2,\Z)$, $c(A)$, and $c(\A)$
					\page9		\endhead
    \head{2.3.} Spines of torus bundles \page{14}	\endhead
    \head{2.4.} Digression: spines of lens spaces
					\page{20}	\endhead
\specialhead{\S 3.} Lower bound for $C^1$-smooth spines transversal to fibers
							\endspecialhead
    \head{3.1.} Morse transformations in simple polyhedra
					\page{23}	\endhead
    \head{3.2.} $\th$-curves in the fibers
					\page{26}	\endhead
    \head{3.3.} Proof of Theorem~12	\page{30}	\endhead
\specialhead{\S 4.} Main conjectures			\endspecialhead
    \head{4.1.}	Triangulations of $T^2$ coming from spines
					\page{31}	\endhead
    \head{4.2.} Can hyperbolic geometry help?		
					\page{37}	\endhead
    \head{4.3.} Other 3-dimensional manifolds
    					\page{40}	\endhead
\specialhead{References}		\page{42}	\endspecialhead
\endtoc
\endtopmatter

\document
\baselineskip12pt minus .1pt

\newpage
\head\S1. Introduction \endhead

\subhead 1.1. Brief outlook \endsubhead

The notion of {\it complexity\/} of three-dimensional manifolds was introduced
by S.~Matveev in~1990, see~\cite{15}. This complexity is a natural
``filtration'' on the set of compact 3-manifolds. It is additive with respect
to taking the connected sum of manifolds, and for any $k\in\Z$ there are only
finitely many compact prime 3-manifolds of complexity at most~$k$; they can be
enumerated by a simple algorithm (however, most of them appear many times in
the list obtained, for example, in the list in~\cite{14,~\S5.2}). For any
compact prime 3-manifold $M$ different from $S^3$, $\R P^3$, $L_{3,1}$, and
$S^2\times S^1$, the complexity $c(M)$ is nothing but the minimal possible
number of tetrahedra in a singular triangulation of~$M$. On the other hand,
these four manifolds are the only closed prime manifolds of complexity~0.

The problem of evaluating the complexity of 3-manifolds is unresolved and
appears to be very difficult. The only manifolds of known complexity are
those with complexity less than or equal to~7. Their lists in~\cite{14}
and~\cite{19} are obtained by enumeration of all special spines of closed
orientable 3-manifolds of small complexity (by the algorithm mentioned above),
followed by determining which of the spines obtained are equivalent (that is,
are spines of the same manifold). This ``equivalence problem'' is difficult.

Obviously, any almost simple spine (or singular triangulation) of a
manifold~$M$ provides an upper bound for~$c(M)$. There is an algorithm for
simplification of a given spine, see~\cite{13}; for all manifolds
from~\cite{14} and~\cite{19}, this algorithm is efficient, that is, stops at a
minimal spine of a manifold. There is no proof of efficiency of this algorithm
in the general case, although one can, of course, use it to find quite
reasonable upper bounds for~$c(M)$.

Much less is known about lower bounds. Clearly, $c(M)>7$ whenever $M$ is not
homeomorphic to any manifold from tables in~\cite{14,~19}. Also, one can easily
show that $c(M)\ge b_1(M,G)-1$ for any commutative group $G$; here $b_1$ is the
first Betti number. However, in most cases these estimates are very inadequate.
Up to now, the only known way to prove that $c(M)=k$ for some $k>0$, where $M$
is a closed prime three-manifold, is to construct a special spine of $M$ with
$k$ vertices (or a singular triangulation of $M$ with $k$ tetrahedra) and
verify that $M$ is not homeomorphic to any manifold of lower complexity.

Here we study 3-manifolds that can be fibered over the circle with torus fiber.
We present ``reasonable'' special spines of these manifolds and discuss
different ways of obtaining lower bounds for their complexity. The conjecture
that arised is very similar to S.~Matveev's conjecture about the complexity of
the lens spaces.

The paper is organized as follows: in the rest of \S1, we give necessary
definitions (following mainly~\cite{14}) and discuss Matveev's conjecture on
the complexity of lens spaces, which implies, by the way, an unexpected theorem
on the running time of the Euclid algorithm (Theorem~3 in~\S1.3). In \S2.3, we
construct pseudominimal spines with small number of vertices for the total
spaces of torus bundles over the circle; as a byproduct, in~\S2.4 we obtain
another description of pseudominimal spines of lens spaces. Both constructions
are based on the study of \thcs in tori, $\SL(2,\Z)$-action on the set $\Sigma$
of their isotopy classes, flips, and the distance $d$ on~$\Sigma$ defined by
flips, see \S2.1 and~\S2.2; another point of view on the flip-distance is
presented in~\S4.2.

We conjecture that the spines constructed in \S2.3 and~\S2.4 are minimal, that
is, the complexity of the corresponding manifolds is equal to the number of
vertices in those spines, and the rest of the paper is devoted to this
conjecture. In \S3 we restrict ourselves to the spines of these spaces that are
transversal to the fibers and prove that the number of vertices of any special
spine with this property is at least one fifth of its conjectured value, see
Theorem~12 (which holds for all Stallings manifolds, not for torus bundles
only, see~\S4.3.2). Till the end of~\S3, everything relies heavily on the study
of $\th$-curves, which form a metric space $(\Sigma,d)$. A generalization of
this construction is presented in~\S4; potentially, it can give a good estimate
or even the exact value of $c(M)$, but requires to solve a problem about
2-chains similar to the problem discussed in~\cite{26,~\S4}; also
see~\cite{10}. Finally, in~\S4.3 we discuss 3-manifolds different from the
total spaces of torus bundles.

Apart from study of $\th$-curves, another approach deserves to be mentioned.
For a 3-manifold $M$, consider its Turaev--Viro invariants, see~\cite{28}. The
examples of these invariants constructed in~\cite{28} involve a root of
unity~$q$, $q^n=1$, of arbitrary degree~$n$. By construction, $TV_q(M)$ is a
certain sum of products over all vertices of~$P$ of some polynomial expressions
in~$q$ assigned to the vertices of~$P$, where $P$ is a special spine of~$M^3$.
Obviously, the degree of $TV_q(M)$ in~$q$ (or in $\sigma=q+q^{-1}$) does not
exceed the maximal degree $d(q)$ (which depends on~$q$ only) of an expression
corresponding to a vertex of~$P$ times the number of vertices of~$P$, where a
spine~$P$ can be assumed to have the minimal possible number of vertices. So it
is possible, in principle, to estimate the complexity of~$M$ by dividing the
degree of $TV_q(M)$ (in~$q$ or in~$\sigma$) over~$d(q)$; since $q$ is an
$n$\snug th root of unity, it makes sense to consider the limit of that ratio
as $n$ tends to infinity. We do not elaborate this idea in the present paper;
note that evaluating the invariants $TV_q(M)$ is a problem difficult enough by
itself, and spines constructed in~\S2.3 can be used in the calculations.

{\bf Acknowledgements.} The idea to relate \thcs and their flips (studied
in~\cite2) to spines of fibered spaces and their complexity is due to
V.~Turaev (see Fig.~10). R.~Fernandes suggested a simplification of the
original proof of Theorem~3. The author would also like to thank
Yu.~Baryshnikov, Yu.~Burman, V.~Nikulin, M.~Polyak, V.~Vassiliev, and
especially S.~Matveev for many useful discussions. The author thanks Instituto
Superior T\'ecnico (Lisbon, Portugal), Isaac Newton Institute (Cambridge,
U.K.), and IHES (Bures-sur-Yvette, France) for their kind hospitality, and the
former Institute for the support of S.~Matveev's visit, which was very
important for this work.

\subhead 1.2. Definitions \endsubhead

In this section we follow~\cite{14,~15}.  By $K$ denote the 1-dimensional
skeleton of the tetrahedron, which is nothing but the clique (that is, the
total graph) with 4 vertices. Note that $K$ is homeomorphic to a circle with
three radii.

\definition{Definition~1} A compact 2-dimensional polyhedron is called {\it
almost simple\/} if the link of each of its points can be embedded in~$K$. An
almost simple polyhedron $P$ is said to be {\it simple\/} if the link of each
point of~$P$ is homeomorphic to either a circle or a circle with a diameter or
the whole graph~$K$. A point of an almost simple polyhedron is {\it
non-singular\/} if its link is homeomorphic to a circle, it is said to be a
{\it triple point\/} if its link is homeomorphic to  a circle with a diameter,
and it is called a {\it vertex\/} if its link is homeomorphic to~$K$. The set
of singular points of a simple polyhedron~$P$ (i.e., the union of the vertices
and the triple lines) is called its {\it singular graph\/} and is denoted
by~$SP$.							\enddefinition

It is easy to see that any compact subpolyhedron of an almost simple polyhedron
is almost simple as well. Neighborhoods of non-singular and triple points of a
simple polyhedron are shown in Fig.~1\,a,\,b; Fig.~1\,c--f represent four
equivalent ways of looking at vertices.

\midinsert

\epsfxsize=300pt

\centerline{\epsffile{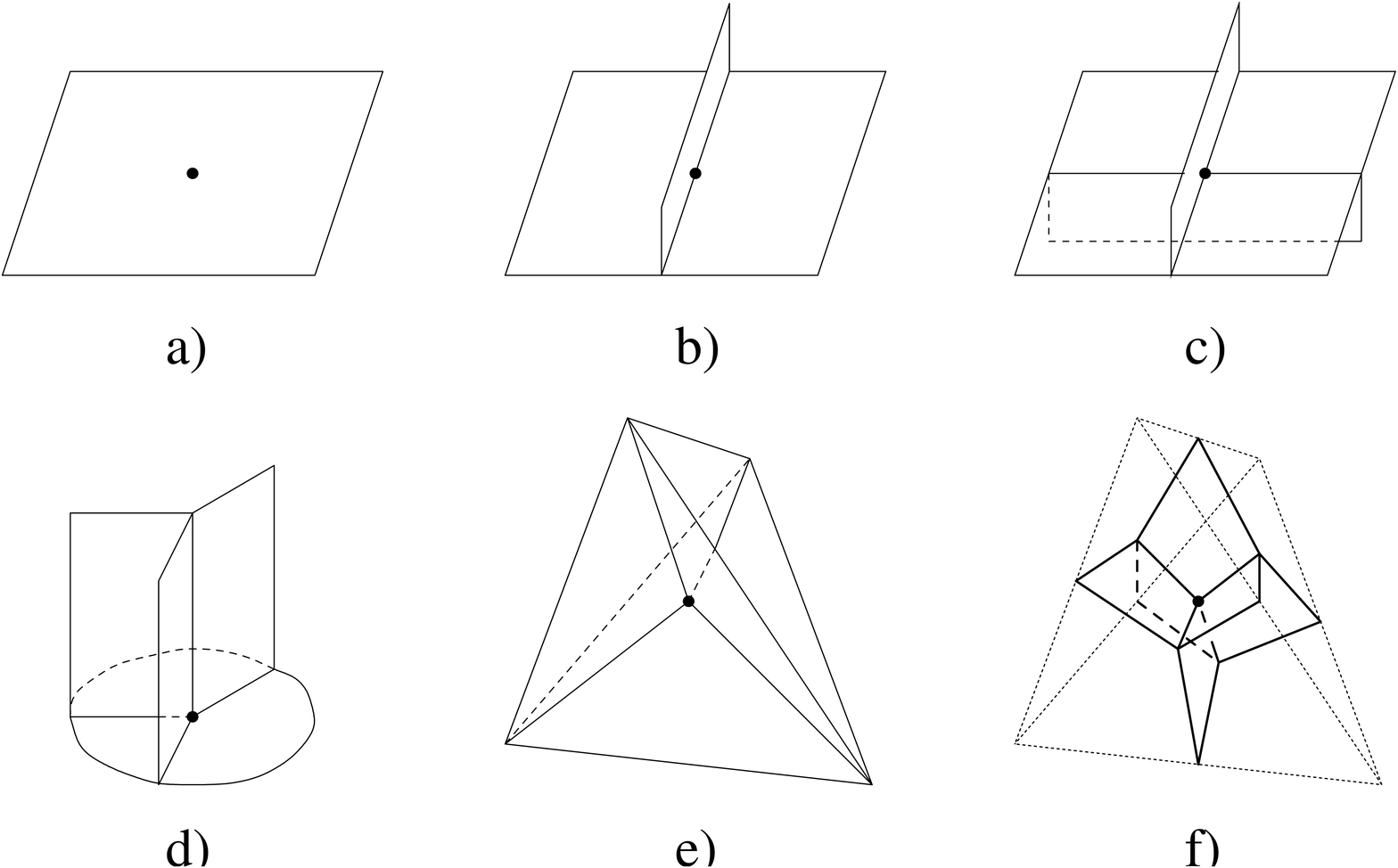}}

\botcaption{Figure 1} Nonsingular (a) and triple (b) points; ways of looking
at vertices (c--f)						\endcaption

\endinsert

\definition{Definition~2} A simple polyhedron~$P$ with at least one vertex is
said to be {\it special\/} if it contains no closed triple lines (wihtout
vertices) and every connected component of $P\sm SP$ is a 2-dimensional cell.
								\enddefinition

\definition{Definition~3} A polyhedron $P\subset\Int M$ is called a {\it
spine\/} of a compact 3-di\-men\-sional manifold~$M$ if $M\sm P$ is
homeomorphic to $\dd M\times(0,1]$ if $\dd M\ne0$ or to an open 3-cell if $\dd
M=0$. In the other words, $P$ is a spine of~$M$ if a manifold $M$ with boundary
(or punctured at one point closed manifold~$M$) can be collapsed onto~$P$. A
spine $P$ of a 3-manifold $M$ is said to be {\it almost simple}, {\it simple},
or {\it special\/} if it is an almost simple, simple, or special polyhedron,
respectively.							\enddefinition

\definition{Definition~4} The {\it complexity\/} $c(M)$ of a complact
3-manifold~$M$ is the minimal possible number of vertices of an almost simple
spine of~$M$. An almost simple spine with the minimal possible number of
vertices is said to be a {\it minimal\/} spine.			\enddefinition

\proclaim{Theorem~1~\cite7} Any compact \rom3-manifold has a special
spine.								\endproclaim

\proclaim{Theorem~2~\cite{15}} Let $M$ be a compact orientable prime
\rom3-manifold with in\-com\-press\-ible \rom(or empty\rom) boundary and
without essential annuli. If $c(M)>0$ \rom(that is, if $M$ is different from
\rom(possibly punctured\rom) $S^3$, $\R P^3$, $L_{3,1}$, and $S^2\times S^
1$\rom), then any minimal almost simple spine of~$M$ is special.  \endproclaim

Recall that a 3-manifold $M$ is called prime if it cannot be represented as a
connected sum $M=M_1\#M_2$ with $M_1$, $M_2$ both different from~$S^3$.

\remark{Remark~1} In this paper, we consider lens spaces $L_{p,q}$, $q>3$, and
the total spaces of torus bundles over the circle. All these manifolds satisfy
the assumptions of Theorem~2.					\endremark

\remark{Remark~2} Starting from a special spine $P$ of a manifold $M$, one can
triangulate~$M$ into $n$ tetrahedra, where $n$ is the number of vertices of~$P$.
This singular triangulation has the only vertex somewhere inside $M\sm P$, its
edges are dual to the 2-cells of~$P$, and triangles are dual to the edges
of~$P$. On the other hand, given a singular triangulation of~$M$ containing $n$
tetrahedra, one can easily obtain a special spine of the manifold~$M$ punctured
at all vertices of the triangulation. It was shown in~\cite{15} that puncturing
does not affect the complexity. Thus for a manifold $M$ satisfying assumptions
of Theorem~2 (in particular, for any prime manifold without boundary), its
complexity $c(M)$ is equal to the minimal possible number of tetrahedra in a
singular triangulation of~$M$, provided that~$c(M)>0$.		\endremark

\remark{Remark~3} Let a special spine $P$ of a manifold $M$ without boundary
have $n$ vertices. Since each vertex of the graph $SP$ has degree~4, $P$
contains $2n$ edges. Since the Euler characteristic of any 3-manifold equals
zero and $M\sm P$ is a 3-cell, we have the equality $n-2n+f-1=0$, which implies
$f=n+1$, where $f$ stands for the number of 2-dimensional ``faces'' of~$P$. It
follows from the construction of Remark~2 that the groups $\pi_1(M)$ and $H_1(
M)$ have at most $f$ generators. Therefore, $f-1=c(M)\ge b_1(M)-1$.  \endremark

\subhead 1.3. Example: lens spaces \endsubhead

\definition{Definition~5} Let $p,q$ be coprime positive integers. The {\it 
Euclid complexity\/} $E(p,q)$ is the number of subtractions (not divisions!)
that the Euclid algorithm takes to convert the pair $(p,q)$ into the pair
$(0,1)$. It is easy to see that $E(p,q)$ equals the sum of the denominators of
the continued fraction representing any of the rational numbers $p/q$
and~$q/p$.							\enddefinition

A good exposition of the Euclid algorithm and continued fraction theory can be
found in~\cite{7,~29}.

\proclaim{Conjecture~1~\cite{14,~15}} The complexity of the lens space $L_{p,q
}$ is equal to $c(L_{p,q})=E(p,q)-3$.				\endproclaim

Pseudominimal special spines of the spaces $L_{p,q}$ with $E(p,q)-3$ vertices
were constructed in~\cite{13,~14}. Pseudominimality of a spine means that no
simplification move can be applied to it; for exact definitions, see~\cite{14}.
In the end of \S2 we present another construction of these spines. Note that
the manifolds $L_{p,q}$ and $L_{p,p-q}$ are homeomorphic, and so are the
manifolds $L_{p,q}$ and $L_{p,r}$, where $0<q,r<p$ and $qr\equiv1\mod p$. So
Conjecture~1 implies that $E(p,q)=E(p,p-q)$ and $E(p,q)=E(p,r)$ for $p,q,r$ as
above; if these corollaries did not hold, Conjecture~1 would fail automatically.
However, they are true. Indeed, $E(p,q)=E(q,p-q)+1$ and $E(p,p-q)=E(p-q,q)+1$,
which implies $E(p,q)=E(p,p-q)$. The second corollary is a true statement, too,
but this is far less obvious.

\proclaim{Theorem~3} Let $0<q,r<p$ and $qr\equiv\pm1\mod p$. Then $E(p,q)=E(p,r
)$.								\endproclaim

\demo{Proof} We can suppose that $p\ge3$. Let us introduce two {\it row
transformation\/} matrices
	$$R_1=\pmatrix1&1\\ 0&1\endpmatrix\qquad\hbox{and}
			\qquad R_2=\pmatrix1&0\\ 1&1\endpmatrix.$$
%
Obviously, we have
	$$R_1^{\pm1}\pmatrix a&b\\ c&d\endpmatrix=
	\pmatrix a\pm c&b\pm d\\ c&d\endpmatrix\qquad\hbox{and}\qquad
	R_2^{\pm1}\pmatrix a&b\\ c&d\endpmatrix=
	\pmatrix a&b\\ a\pm c&b\pm d\endpmatrix.$$
Consider the expansion of $p/q$ in a continued fraction
	$$\frac pq=n_1+\cfrac1\\ n_2+\cfrac1\\ \ddots\lower6pt
	\hbox{$+\dfrac1{n_k}$}\endcfrac\,.$$
Set $U=R_2^{-1}R_\e^{-n_k+1}\dots R_1^{-n_3}R_2^{-n_2}R_1^{-n_1}$, where $\e=1$
for $k$ odd and $\e=2$ for $k$ even; note that $n_1\ge1$ and $n_k\ge2$. It is
easy to see that $U$ takes vector $(p,q)^{\T}$ to $(1,0)^{\T}$ (where $\T$
stands for transposing): $R_1^{-n_1}$ takes $(p,q)^{\T}$ to $(p-n_1q,q)^{\T}$
and so forth, according to the Euclid algorithm, the only exception being that
at the last step we apply $R_2^{-1}$ to $(1,1)^{\T}$, not $R_1^{-1}$, in order
to convert $(1,1)^{\T}$ to $(1,0)^{\T}$, not to $(0,1)^{\T}$.

First suppose that $qr\equiv-1\mod p$, thus $qr=sp-1$ for some positive
integer~$s$. Since $1\le q<p$ and $1\le r<p$, we have $s\le q$ and $s\le r$.
Let us consider the inverse matrix
	$$U^{-1}=R_1^{n_1}R_2^{n_2}R_1^{n_3}\dots R_\e^{n_k-1}R_2.\tag1$$
We proclaim that
	$$U^{-1}=\pmatrix p&r\\ q&s\endpmatrix.$$
Indeed, equation~(1) implies that $U^{-1}$ has the following properties:
\roster
\item"1)" the first column of $U^{-1}$ is $(p,q)^{\T}$;
\item"2)" the determinant of $U^{-1}$ equals~1;
\item"3)" the second column entries of $U^{-1}$ are positive, and 
\item"4)" they do not exceed the corresponding first column entries.
\endroster
Property 1 follows from the construction of the sequence involved in~(1), and
property 2 is obvious. It is clear that both second column entries of $U^{-1}$
are nonnegative; they are both positive, because the right hand side of~(1)
contains both $R_1$ and~$R_2$. The last property holds for~$R_2$ and survives
under left multiplications by~$R_1$ and~$R_2$. The first two properties imply
that the second column of $U^{-1}$ is $(r+mp,s+mq)^{\T}$ for some $s\in\Z$,
and the last two properties show that in fact $m=0$.

Note that $E(p,q)=n_1+n_2+\ldots+n_k$ equals the sum of the exponents in~(1).
Since $R_1^{\T}=R_2$ and $R_2^{\T}=R_1$, for the transposed inverse matrix
$(U^{-1})^{\T}$ we have
	$$\pmatrix p&q\\ r&s\endpmatrix=(U^{-1})^{\T}=R_1R_{3-\e}^{n_k-1}
			\dots R_2^{n_3}R_1^{n_2}R_2^{n_1},\tag2$$
where $n_k\ge2$ and $n_1\ge1$. This yields the sequence of $n_1+n_2+\ldots+n_k=
E(p,q)$ subtractions that converts the pair $(p,r)$ into the pair $(1,0)$. Such
a sequence is unique up to a possible application of several subtractions $R_1$
to the column $(1,0)^{\T}$ (clearly, $R_1$ does not change it); it is nothing
but the Euclid algorithm provided that there are no those ``fake'' subtractions.
This condition is satisfied, because the last subtraction (the inversed
rightmost one in~(2), or the leftmost one in a similar expression for $U^{\T}$)
is $R_2$, not~$R_1$. Therefore $E(p,q)=E(p,r)$ for $qr\equiv-1\mod p$. In the
case $qr\equiv1\mod p$, note that $q(p-r)\equiv-1\mod p$ and $0<p-r<p$.
Consequently, $E(p,q)=E(p,p-r)=E(p,r)$.				\qed\enddemo

Theorem~3 means that Conjecture~1 passes a nontrivial ``sanity test''.

\head\S2. Upper bound \endhead

\subhead 2.1. $\th$-curves	\endsubhead

Let us recall some definitions and results from~\cite2\footnote{We do not
follow the notation of~\cite2 here.}.

\definition{Definition~6} A {\it $\th$-curve\/} $L\subset T^2$ is a graph with
two vertices and three edges (not loops) connecting these vertices, embedded
in~$T^2$ in such a way that the edges are pairwise non-homotopic; this is
equivalent to the condition that the complement $T^2\sm L$ is a 2-dimensional
cell.								\enddefinition

\midinsert

\epsfxsize=120pt

\centerline{\epsffile{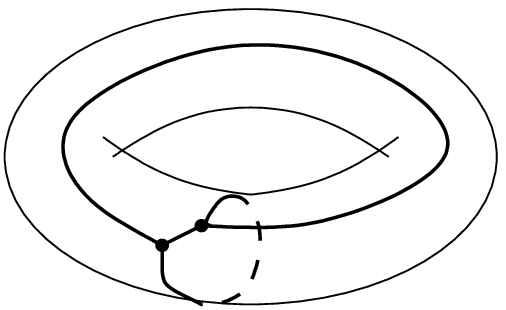}}

\botcaption{Figure 2} A $\th$-curve	\endcaption

\endinsert

Up to an isotopy, any two \thcs can be taken to one another by a linear
automorphism of the torus, see~\cite2. Another way to change the isotopy class
of a \thc is to apply a sequence of flips.

\definition{Definition~7} A {\it flip} along an edge of a trivalent graph (in
particular, of a $\th$-curve) is an invertible restructuring of the graph that
acts on a neighborhood of this edge as shown on~Fig.~3. A flip does not change
the graph outside of this neighborhood.				\enddefinition

\midinsert

\epsfxsize=250pt

\centerline{\epsffile{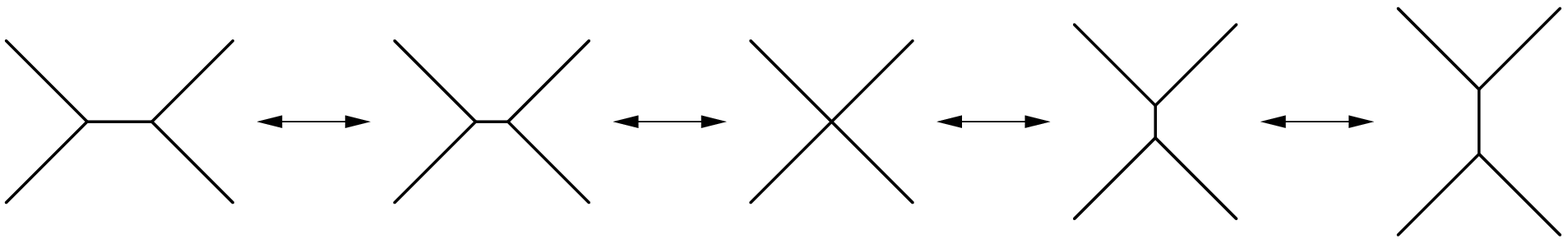}}

\botcaption{Figure 3} A flip	\endcaption

\endinsert

For any two \thcs $L_1$, $L_2$, there exists a sequence of flips (and
isotopies) that takes $L_1$ to~$L_2$, see~\cite{1,~2}. Now let us recall how
one can find the minimal number of flips required for such a sequence.

For a \thc $L$, there are three unoriented (or six oriented) cycles in $\pi_1
(T^2)$ formed by pairs of edges of~$L$. These cycles also can be represented by
the three 1-cells of the singular triangulation of~$T^2$ dual to the cell
decomposition defined by~$L$. The six points in the lattice $\Z^2=\pi_1(T^2)=
H_1(T^2,\Z)$ corresponding to these cycles are the vertices of some convex
centrally symmetric hexagon~$W(L)$.

\definition{Definition~8} The hexagon $W(L)$ is said to be {\it associated\/}
to a $\th$-curve~$L$. A hexagon $W$ with the vertices in $\Z^2$ is called
{\it admissible\/} if it is associated to some $\th$-curve. The {\it standard
hexagon\/} is the hexagon $W_0$ with the vertices $\pm(1,0)$, $\pm(0,1)$, and
$\pm(1,-1)$, see Fig.~4. It is associated to the \thc shown on Fig.~2 (under
a natural choice of a parallel and a meridian of $T^2$ as a basis of~$H_1(T^2)
=Z^2$).								\enddefinition

\proclaim{Theorem~4~\cite2}

\rom{1)} A hexagon $W$ \rom(with vertices at lattice points\rom), centrally
symmetric with respect to the origin~$O$, is admissible if and only if it has
the following properties\rom:
\roster
\item"a)" if $X$ and $Y$ are nonopposite vertices of~$W$, then the area of the
	triangle $OXY$ equals~$1/2$\rom;
\item"b)" if $X$, $Y$, and $Z$ are three consecutive vertices of~$W$, then
	$\oa{OY}=\oa{OX}+\oa{OZ}$.
\endroster
Properties \rom{a)} and \rom{b)} are equivalent. The origin is the only
interior lattice point of an admissible hexagon~$W$. The vertices of~$W$ are
the only lattice points on its boundary.

\rom{2)} Two \thcs $L$ and $L'$ are isotopic if and only if $W(L)=W(L')$. For
any two \thcs $L$ and $L'$ there exists an operator $A\in\SL(2,\Z)$ such that
$AL$ is isotopic to~$L'$ and $A(W(L))=W(L')$.

\rom{3)} Let \thcs $L$ and $L'$ differ by the flip along an edge~$e$. Then
associated hexagons $W$ and $W'$ have in common two pairs of opposite
vertices that correspond to cycles $\sigma$ and~$\mu$ dual to two other edges.
The remaining pair of the vertices is $\pm(\sigma+\mu)$ for one of the hexagons
$W$, $W'$ and $\pm(\sigma-\mu)$ for the other one.		\endproclaim

Since a \thc has three edges, three different flips can be applied. Fig.~4
shows how they change the standard hexagon~$W_0$. The result of a flip
transformation of an arbitrary admissible hexagon can be represented by the
same picture with another coordinate system, because any admissible hexagon is
$\SL(2,\Z)$-equivalent to~$W_0$.

\midinsert

\epsfxsize=250pt

\centerline{\epsffile{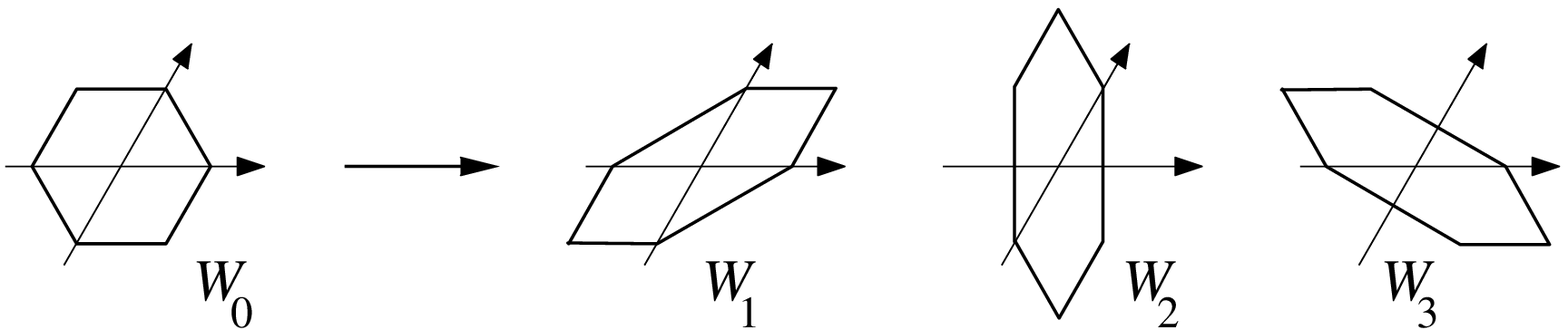}}

\botcaption{Figure 4} Action of three flips on the standard hexagon \endcaption

\endinsert

According to the second part of Theorem~4, we can study sequences of flips
that convert an admissible hexagon~$W_1$ into another admissible hexagon~$W_2$
rather than sequences of flips that take a \thc $L_1$ into another
$\th$-curve~$L_2$. So it is natural to introduce a graph~$\Gamma$ that has the
admissible hexagons as its vertices and flips as its edges; the number of flips
required to convert $W_1$ into $W_2$ equals the distance between the
corresponding vertices of~$\Gamma$. Clearly, $\Gamma$ is a trivalent graph. It
turns out that $\Gamma$ is a tree, see~\cite2. More information about~$\Gamma$
can be found at~\cite{25, Ch.~II,~\S1}.

\definition{Definition~9} A {\it leading vertex\/} of an admissible hexagon~$W$
is its vertex that is the most distant from the origin with respect to the
quadratic form $Q(x,y)=x^2+xy+y^2$.				\enddefinition

The standard hexagon~$W_0$ is a unit regular hexagon with respect to $Q(x,y)$.
Any other admissible hexagon has only one pair of opposite leading vertices.

\proclaim{Theorem~5~\cite2} Let $(p,q)$ be a leading vertex of an admissible
hexagon~$W\ne W_0$. Suppose that $p>0$ and $q>0$. Then $d(W,W_0)=E(p,q)$, where
$d(W,W_0)$ stands for the distance between $W$ and $W_0$ in~$\Gamma$ and
$E(p,q)$ is the Euclid complexity, see Definition~\rom5 above.	\endproclaim

In fact, the steps (flips) of the only way from (the vertex of~$\Gamma$
corresponding to) $W_0$ to (the vertex corresponding to) $W$ in~$\Gamma$ are in
a natural one-to-one correspondence with the steps (subtractions) of the Euclid
algorithm applied to the pair $(p,q)$. There is an algorithm that constructs
the path from $W$ to~$W_0$: start at~$W$ and apply the flip that decreases
the length of a hexagon; such a flip is unique unless the hexagon is~$W_0$. The
numbers $p$, $q$ are coprime by virtue of the first part of Theorem~4, and for
any pair of coprime numbers $(p,q)$, there exists exactly one admissible
hexagon with a leading vertex at $(p,q)$, cf.~\cite2. For the detailed proof of
Theorem~5, see~\cite2.

If the coordinates $(p,q)$ of a leading vertex of~$W$ are both negative, one
should consider the opposite vertex, or rotate the coordinate system by~$\pi$.
If $pq<0$, one should rotate the coordinate system by $\pm\pi/3$ (we consider
``triangular'' coordinates shown on Fig.~4 rather than rectangular
coordinates), that is, to apply the coordinate change $\pmatrix 0&-1\\ 1&\ph-1
\endpmatrix$ or its inverse, to make both coordinates of one of the
leading vertices positive. Then Theorem~5 may be applied. This gives the
following simple answer: $d(W,W_0)=E(|p|,|q|)-1$ whenever $pq<0$, where $(p,q)$
is a leading vertex of~$W$. The $-1$ summand can be explained as follows: if,
for example, $p<0$ and $q>0$, the process of converting $W$ into $W_0$ by flips
corresponds to the converting of the unordered pair $(-p,q)$ to the pair
$(1,1)$, not to $(0,1)$, by subtractions according to the Euclid algorithm
(since we can consider $(-1,1)$ as a leading vertex of $W_0$); of course, this
takes one subtraction less. Also see~\S4.2.

\definition{Definition~10} For a matrix $A\in\SL(2,\Z)$ we define its {\it
complexity\/} $c(A)$ by putting $c(A)=d(W_0,AW_0)$.	\enddefinition

To calculate the number $c(A)$, find a leading vertex $(p,q)$ of the hexagon~$A
W_0$. Then we get $c(A)=E(|p|,|q|)$ if $pq>0$ or $c(A)=E(|p|,|q|)-1$ if $pq<0$.
In particular, if $AW_0=W_0$ (there are six matrices $A$ with this property),
we get $c(A)=0$.

\subhead 2.2. On $\SL(2,\Z)$, $c(A)$, and $c(\A)$ \endsubhead

Since all admissible hexagons are $\SL(2,\Z)$-equivalent and the action of this
group preserves the existence of a flip connecting two given hexagons (thus,
there is an action of $\SL(2,\Z)$ on the graph~$\Gamma$), we have $d(W_1,W_2)=
d(BW_1,BW_2)$ for $B\in\SL(2,\Z)$. In particular, 
  $$c(A^{-1})=d(W_0,A^{-1}W_0)=d(AW_0,AA^{-1}W_0)=d(W_0,AW_0)=c(A)\tag3$$
and $c(A)=d(W_0,AW_0)=d(BW_0,BAW_0)$, which may be different from $d(BW_0,\amb
ABW_0)=d(W,AW)$. Thus we have $c(A)=d(W,BAB^{-1}W)$ for any admissible
hexagon~$W=BW_0$.

It is not true that $c(A)=d(W,AW)$ for any admissible hexagon~$W$. This means
that the number $c(A)$ is not a conjugacy class invariant, contrary to a
statement contained implicitly in~\cite2.

\example{Example} Let $A=\pmatrix\ph-171&\ph-100\\ -289&-169\endpmatrix$ and $B
=\pmatrix\ph-10&-17\\ -17&\ph-29\endpmatrix$. Note that $A,B\in\SL(2,\Z)$. By a
straightforward calculation, we obtain $c(A)=13$, while for a conjugate matrix
$A'=B^{-1}AB=\pmatrix1&1\\ 0&1\endpmatrix$ we have $c(A')=1$.	\endexample

The following problem arises: find the minimal value of $c(A)$ over the whole
conjugacy class of~$A$ in~$\SL(2,\Z)$.

\definition{Definition~11} The {\it complexity\/} of an operator $\A\in\SL(2,\Z
)$ is the minimal possible complexity of matrices that represent~$\A$ in all
possible bases of the lattice~$\Z^2$:
$$c(\A)=\min_{A\sim A_0}c(A)=\min_{B\in\SL(2,\Z)}c(B^{-1}A_0B),$$
where $A_0$ is any matrix representing $\A$ and $\sim$ denotes conjugacy
in~$\SL(2,\Z)$. In other words, $c(\A)=\min d(W,\A W)$ over all admissible
hexagons~$W$. A matrix $A$ of an operator $\A$ in some basis is said to be a
{\it minimal\/} matrix of~$\A$ if $c(A)=c(\A)$, that is, if $c(A)\le c(A')$ for
any matrix $A'$ conjugated to~$A$. An admissible hexagon $W$ is called {\it
minimal\/} for~$\A$ if $d(W,\A W)=c(\A)$.			\enddefinition

In~\S3, we will be interested in the sequence $\{c(\A^k)\}$. Properties of
this sequence depend on the trace of~$\A$. Recall (see~\cite{6,~\S0}) that the
operator $\A$ is called {\it elliptic\/} if $|\tr\A|<2$. In this case $\tr\A=0$
or $\tr\A=\pm1$, and the equation $\A^2-(\tr\A)\A+(\det\A)I=0$ implies either
$\A^2=-I$ or $\A^3\pm I=0$, because $\det\A=1$. Thus elliptic operators are
periodic of period 3, 4 or~6, and so are the sequences $\{c(\A^k)\}$; both
eigenvalues of an elliptic operator are roots of unity. If $\tr\A=\pm2$, we say
that $\A$ is a {\it parabolic\/} operator. In this case $(\A\pm I)^2=0$ and
$\A$ is $\SL(2,\Z)$ conjugated to either $\pmatrix 1&n\\ 0&1\endpmatrix$ or
$\pmatrix-1&n\\ \ph-0&-1\endpmatrix$, where $n\in\Z$. So $\A$ is either a
periodic operator (if $n=0$) or, up to a sign, a power of the Jordan block
$\pmatrix 1&1\\ 0&1\endpmatrix$; both eigenvalues of a parabolic operator equal
$\pm1$. Finally, $\A$ is {\it hyperbolic\/} if $|\tr\A|>2$. In this case the
eigenvalues of $\A$ are different real numbers, and $\A$ is a hyperbolic
rotation. Also see~\cite{22,~\S5}.

\proclaim{Lemma~1} $c(AB)\le c(A)+c(B)$. Moreover, $c(AB)\equiv c(A)+c(B)
\mod2$.								\endproclaim

\demo{Proof} By definition, $c(A)=d(W_0,AW_0)$, $c(AB)=d(W_0,ABW_0)$, and $c(B)
=d(W_0,BW_0)=d(AW_0,ABW_0)$. Since $d$ is a metric on~$\Gamma$, the triangle
inequality $d(W_0,ABW_0)\le d(W_0,AW_0)+d(AW_0,ABW_0)$ holds. Since the graph
$\Gamma$ is a tree, we have $d(W_0,ABW_0)=d(W_0,AW_0)+d(AW_0,ABW_0)-2k$, where
$k\in\Z_{\ge0}$, see Fig.~5.					\qed\enddemo

\midinsert

\epsfxsize=150pt

\centerline{\epsffile{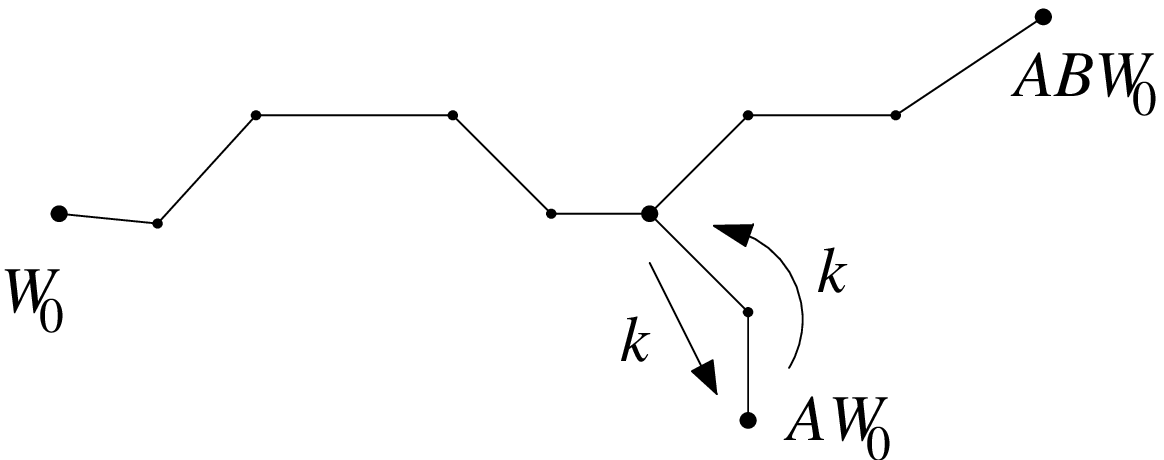}}

\botcaption{Figure 5} Additivity of tree distance parities
\endcaption

\endinsert

\proclaim{Theorem~6} Reduction of $c(A)$ modulo \rom2 coincides with the unique
epimorphism of the modular group $\SL(2,\Z)/\{\pm I\}$ to $\Z_2$.  \endproclaim

\demo{Proof} Lemma~1 implies that this reduction is a homorphism of $\SL(2,\Z)$
to~$\Z_2$. It is an epimorphism since it takes $\pmatrix 1&1\\ 0&1\endpmatrix$
to~$1\mod2$. It is an epimorphism of the modular group $G=\SL(2,\Z)/\{\pm I\}$,
because $c(A)=c(-A)$: the hexagon $AW_0$ is centrally symmetric, whence $AW_0=
-AW_0$. Such an epimorphism $\f$ is unique. Indeed, it is defined by its values
$\f(S)$ and $\f(T)$ on the elements $S=\pmatrix 0&-1\\ 1&\ph-0\endpmatrix$ and
$T=\pmatrix1&1\\ 0&1\endpmatrix$, which generate the modular group,
see~\cite{24, Ch.~VII,~\S1}. The relation $(ST)^3=1$ in $G$ implies that $\f(S)
=\f(T)$ in~$Z_2$. Since $\f$ is not identically zero, we have $\f(S)=\f(T)=1$.
								\qed\enddemo

Let us explain how to find a minimal matrix of an operator~$\A$. Let $A$ be its
matrix in some basis, $W_0$ be the standard hexagon in this basis, $W=AW_0$,
and $W_0,W_1,\dots,W_{c(A)}=W$ be the shortest path from $W_0$ to $W$
in~$\Gamma$. Then $W=AW_0,AW_1,\dots,AW_{c(A)}=AW=A^2W_0$ is the shortest path
from $W$ to $AW$ in~$\Gamma$. Both $W_{c(A)-1}$ and $AW_1$ are neighbors of
a trivalent vertex $W$ of~$\Gamma$. Compare them.

\proclaim{Theorem~7} Any matrix $A$ with $c(A)\le1$ is minimal. A matrix $A$
with $c(A)>1$ is minimal if and only if $W_{c(A)-1}\ne AW_1$.	\endproclaim

\demo{Proof} If $c(A)=0$, the matrix $A$ is minimal. Let $c(A)=1$, whence $c(\A
)\le1$. It follows from eq.~(3) and Lemma~1 that $c(B^{-1}AB)\equiv c(A)+2c(B)
\equiv c(A)=1\pmod2$. So $c(\A)$ is odd, thus $c(\A)\ne0$ and $A$ is minimal.

Suppose that $W_{c(A)-1}=AW_1$. The operator $\A$ takes $W_1$ to $AW_1=W_{c(A)
-1}$. We have $c(\A)\le d(W_1,\A W_1)=d(W_1,W_{c(A)-1})=c(A)-2$ (unless $c(A)
\le1$), i.\,e., the matrix $A$ is not minimal. This proves the ``only if'' part
of the Theorem.

Suppose that the matrix $A$ is not minimal. This implies that the standard
hexagon~$W_0$ is not minimal. Let $V$ be a minimal hexagon for the
operator~$\A$. By $\gamma_0$ denote the path from $V$ to $AV$ in $\Gamma$ of
length~$c(\A)$. For any $k\in\Z$ let $\gamma_k=A^k\gamma$ be the path from
$A^kV$ to $A^{k+1}V$ in $\Gamma$ (recall that $\Gamma$ carries an action of
$\SL(2,\Z)$). Put $\gamma=\bigcup\limits_{k\in\Z}\gamma_k$. Note that any
vertex of $\Gamma$ on $\gamma$ represents a minimal hexagon, so $W_0\in\Gamma
\sm\gamma$. We have to consider three cases.

Case~1: $c(\A)>1$. Then two vertices of $\gamma$ neighboring with $A^kV$, $k\in
\Z$, are different because of the ``only if'' statement proven above. So
$\gamma$ is homeomorphic to a line, because two neighbors of any interior
vertex of any $\gamma_k$ are different, too (since $\gamma_k$ is the shortest
path from $A^kV$ to $A^{k+1}V$), and the graph $\Gamma$ is a tree.

By $W_0U_0$ denote the shortest path from $W_0$ to~$\gamma$ (that is, $U_0$ is
the first point of $\gamma$ belonging to any path from $W_0$ to a point
of~$\gamma$). Then $W\,U$ and $AW\,AU$ are the shortest paths from $W$ and $AW$
to~$\gamma$, where $W=AW_0$ and $U=AU_0$. Obviously, these paths end at
different points of~$\gamma$: if $U_0\in\gamma_k$, then $U\in\gamma_{k+1}$ and
$AU\in\gamma_{k+2}$. Since $\Gamma$ is a tree, the paths $W_0\,U_0$, $W\,U$,
and $AW\,AU$ are mutually disjoint. This means that $W_0\,U_0\,U\,W$ is the
shortest path from $W_0$ to $W$ and $W\,U\,AU\,AW$ is the shortest path from
$W$ to $AW$, see Fig.~6. The leg $WU=A(W_0U_0)$ of this path is not empty.
Consequently, the penultimate vertex of the path $W_0W$ coincides with the
first (different from $W$) vertex of the path $W\,AW$, which proves the
statement of the ``if'' part of the Theorem.

\midinsert

\epsfxsize=200pt

\centerline{\epsffile{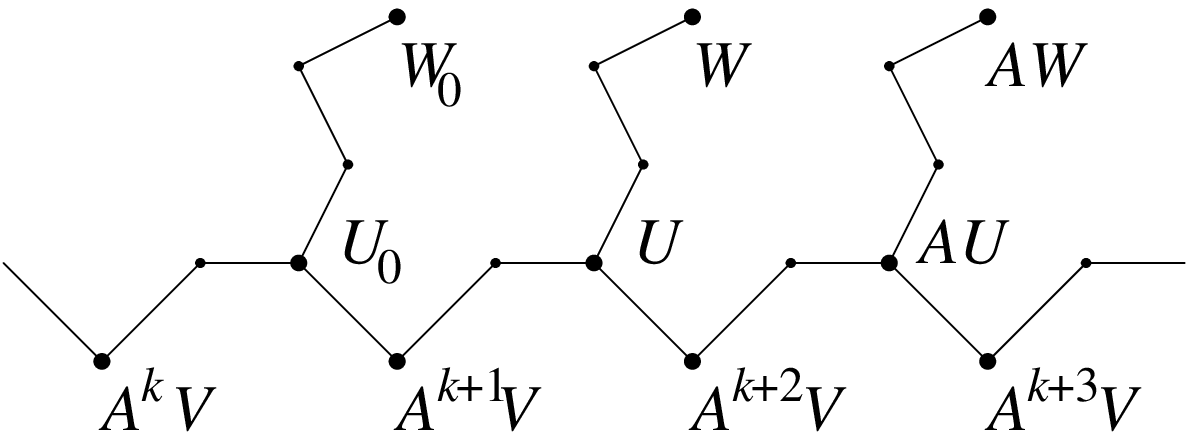}}

\botcaption{Figure 6} Paths $W_0W$ and $W\,AW$ overlap
\endcaption

\endinsert

Case~2: $c(\A)=1$. Let $V$ be a minimal hexagon for~$\A$. If $A^2V\ne V$, then
$A^{k+2}V\ne A^kV$ for any $k\in\Z$, all the $\gamma_k$ are different,
$\gamma$~is homeomorphic to a line, and we can repeat the argument of Case~1.
If $A^2V=V$, we have $A^kV=V$ for $k$ even and $A^kV=AV$ for $k$~odd. Without
loss of generality, it can be assumed that the path from $W_0$ to $V$ does not
pass through $AV$. The transformation $\A$ takes this path to the path from $W=
AW_0$ to $AV$, which does not pass through~$V$. So the shortest path from $W_0$
to~$W$ is $W_0\,V\,AV\,W$, and thus it contains the edge $V\,AV$. Similarly, the
path from $W$ to $AW$ also contains that edge. Therefore, these two paths
overlap, which proves the Theorem in the case~$c(\A)=1$.

Case~3: $c(\A)=0$. There exists an admissible hexagon $V$ such that $AV=V$, but
$AW_0\ne W_0$ because $A$ is not a minimal matrix. Consider paths $W_0V$ and $W
V$ in~$\Gamma$, where $W=AW_0$. Let $U$ be the first (most distant from~$V$)
common point of this paths. Recall that $A$ acts on $\Gamma$ and takes the path
$W_0V$ to $WV$. Thus $AU=U$, $W_0\ne U$, and $A$ takes the path $W_0U$ to~$WU$
and the path $WU$ to~$AW\,U$. So the paths $W_0W=W_0UW$ and $W\,AW=WU\,AW$
overlap over the leg~$WU$.					\qed\enddemo

Now we can present an algorithm that finds the number $c(\A)$ and a minimal
matrix $A$ of an operator~$\A$. Start with any matrix $A$ representing this
operator. Apply the criterion of Theorem~7. If either $c(A)\le1$ or $W_{c(A)-1}
\ne AW_1$, the matrix $A$ is minimal. Otherwise, let $V$ be the last common
vertex of the paths $W\,W_0$ and $W\,AW$. Then $V$ lies on $\gamma$ and is a
minimal hexagon for~$\A$. Choose a basis so that $V$ is the standard hexagon.
The matrix of~$\A$ in this basis is minimal, and $c(\A)$ is equal to complexity
of this matrix.

\proclaim{Corollary} The subgraph $\gamma\subset\Gamma$ constructed in the
proof of Theorem~\rom7 is a line if and only if the operator~$\A$ is not
periodic.						\qed\endproclaim

This condition holds if and only if either $c(\A)\ge2$ or $c(\A)=1$ and $\A^2
\ne-I$. The line~$\gamma$ is the ``mainstream'' of the action of~$\A$
on~$\Gamma$. Minimal hexagons for~$\A$ are exactly those corresponding to the
vertices of the subgraph $\gamma\subset\Gamma$. If $\gamma$ is a line, there
are at most $3\,c(\A)$ different minimal matrices for~$\A$, because any minimal
hexagon yields 6 different bases ($OX_1X_2,\,OX_2X_3,\dots,OX_6X_1$, where $X_1
,\dots,X_6$ are the vertices of a hexagon and $O$ is the origin), bases that
differ by a central symmetry give the same matrix expression of~$\A$, and
hexagons $V$ and $AV$ lead to the same set of matrices of~$\A$.

If $c(\A)=1$, then the minimal matrix for~$\A$ is either one of Jordan blocks
$\pmatrix\pm1&\ph-1\\ \ph-0&\pm1\endpmatrix$, $\pmatrix\pm1&-1\\ \ph-0&\pm1
\endpmatrix$ or the $\pm\pi/2$ rotation matrix $\pmatrix\ph-0&\mp1\\ \pm1&\ph-0
\endpmatrix$. These six matrices belong to six different conjugacy classes
in~$\SL(2,\Z)$. In the case of a Jordan block, there are three minimal matrices
for~$\A$ and an infinite number of minimal hexagons (which lie on the
line~$\gamma$). For a rotation, there are three different minimal matrices
(namely, $\pmatrix-1&-2\\ \ph-1&\ph-1\endpmatrix$, $\pmatrix-1&-1\\ \ph-2&\ph-1
\endpmatrix$, and $\pmatrix0&-1\\1&\ph-0\endpmatrix$ in the case of
counterclockwise rotation) and only two minimal hexagons, which have in common
the four lattice points where a positively definite integer quadratic form $Q_
{\A}(\oa v)=\det(\oa v,\A\oa v)$ (for the clockwise rotation, we set $Q_{\A}(
\oa v)=-\det(\oa v,\A\oa v)$) attains its minimal positive value~1.

If $c(\A)=0$ and $\A\ne\pm I$, the minimal hexagon is unique; its six vertices
are the six lattice points where a positively definite integer quadratic form
$Q_{\A}(\oa v)=\pm\det(\oa v,\A\oa v)$ attains value~1. The minimal matrix
for~$\A$ is also unique. Finally, for $\A=\pm I$, any admissible hexagon is
minimal, while the only minimal matrix is, of course,~$\pm\pmatrix1&0\\0&1
\endpmatrix$. We omit the proofs of the statements of this paragraph and two
preceding ones; most of them are straightforward.

\proclaim{Theorem~8} Let $\A$ be a non-periodic operator. Then\rom:
\roster
\item"\rom{1)}" for any integer $k\ne0$, $A^k$ is a minimal matrix for $\A^k$
	if and only if $A$ is a	minimal matrix for $\A$\rom;
\item"\rom{2)}" $c(\A^k)=|k|c(\A)$ for any $k\in\Z$\rom;
\item"\rom{3)}" for any integer $k\ne0$, we have $c(A^k)=|k|c(\A)+b$, where $b=
	c(A)-c(\A)$ is a nonnegative even number.
\endroster
\endproclaim

\demo{Proof} It follows from the proof of Theorem~7 and Corollary that the path
from $W_0$ to $A^kW_0$ consists of three legs $W_0U_0$, $U_0\,A^kU_0$, and $A^k
U_0\,A^kW_0=A^k(U_0W_0)$. If the first leg (and, simultaneously, the last one)
is empty (contains no edges), matrices $A$ and $\A^k$ are both minimal.
Otherwise, neither $A$ nor $A^k$ is minimal. This proves the first statement
and shows that the mainstreams of $\A$ and $\A^k$ coincide: $\gamma(\A^k)=
\gamma(\A)$. The second statement of the Theorem follows from the first one
whenever $k\ne0$; the case $k=0$ is trivial. The third statement follows from
the description of the path from $W_0$ to $A^kW_0$ in~$\Gamma$, see above.
								\qed\enddemo

We conclude the section on $\SL(2,\Z)$, $c(A)$, and $c(\A)$ with one more way
of looking at the mainstream~$\gamma(\A)$. Let $\A\in\SL(2,\Z)$ be a hyperbolic
rotation, so $|\tr\A|>2$ and eigenvalues of $\A$ are different real numbers.
Since admissible hexagons are centrally symmetric, we may assume that the
eigenvalues of $\A$ are positive. The eigenvalues and the slopes of
eigenvectors are quadratic irrationals, because the discriminant of the
characteristic equation $\lambda^2-(\tr\A)\lambda+\det\A=0$ equals $(\tr\A)^2-4$
and is not equal to a square of an integer whenever $|\tr\A|>2$. Draw through
the origin two lines $l_1$, $l_2$ parallel to eigenvectors. They divide the
plane into four parts. For each of these parts, consider the convex hull of the
lattice points inside it. Since $\A$ preserves $l_i$, it preserves the convex
hulls $h_1,\dots,h_4$, as well as their boundaries, which are infinite
sequences of segments. The group of the integers acts on $\dd h_i$ by taking $x
\in\dd h_i$ to $\A^kx\in\dd h_i$ for $k\in\Z$. An admissible hexagon is minimal
for $\A$ (i.e., belongs to $\gamma(\A)$) if and only if its leading vertex (and
hence all its other vertices) lies on some $\dd h_i$. The path from a minimal
hexagon $W\in\gamma$ to its image $\A W$ corresponds to the period of the
continued fraction expansion of the slope $\a$ of $l_i$, which is a quadratic
irrational number. An $\SL(2,\Z)$ coordinate change does not affect this
period: it takes $\a$ to $\dfrac{a\a+b}{c\a+d}$, where $\pmatrix a&b\\ c&d
\endpmatrix\in\SL(2,\Z)$, changing the beginning of the continued fraction
expansion of a quadratic irrational number without affecting its periodic part.
We leave the proofs of these statements to the reader.

\midinsert

\epsfxsize=240pt

\centerline{\epsffile{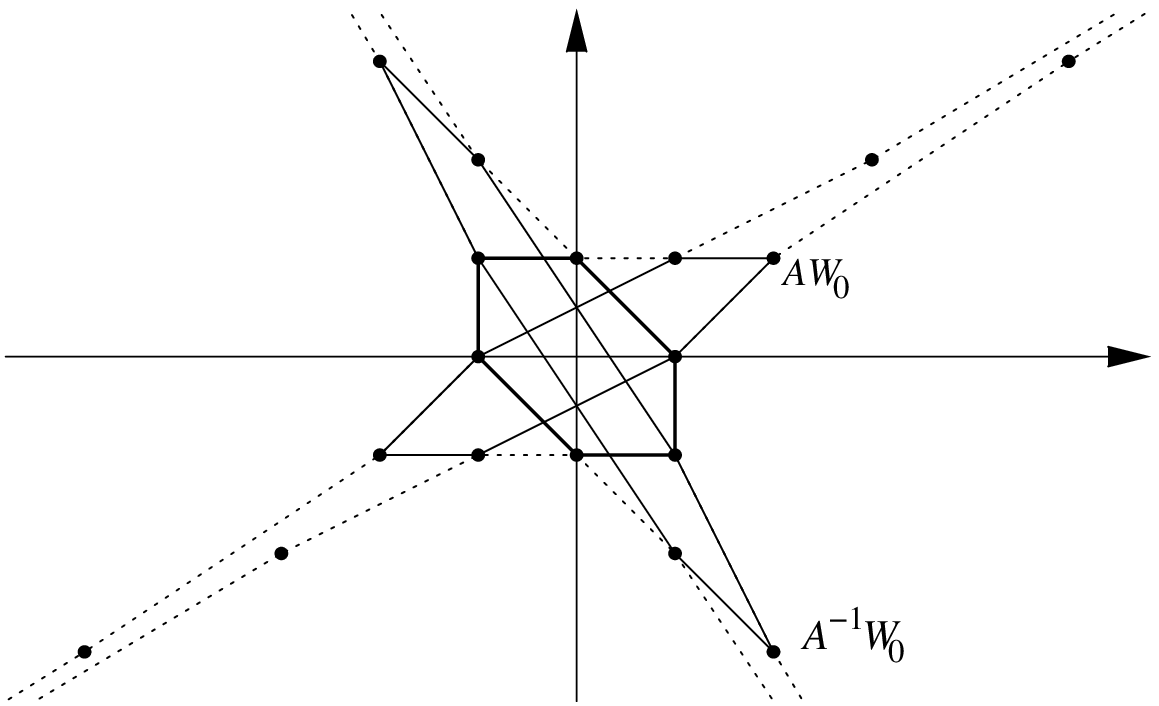}}

\botcaption{Figure 7} The mainstream for a hyperbolic rotation	\endcaption

\endinsert

\example{Example} Let $A=\pmatrix2&1\\1&1\endpmatrix$. This is a minimal matrix
of a hyperbolic operator of complexity~2. The boundaries of the convex hulls~$h
_i$ are represented on Fig.~7 by dotted lines. The coordinates of their corners
in this case are, up to signs, the pairs of consecutive Fibonacci numbers. The
standard hexagon $W_0$, drawn in a bold line, is a minimal hexagon for~$\A$.
The hexagons $A^{-1}W_0$ and $AW_0$ also belong to the mainstream~$\gamma(\A)$.
Since $c(\A)=2$, there are two orbits of the action of the group $\Z$ on
$\gamma(\A)$ defined by the rule $k(W)=A^kW$ for $k\in\Z$ and $W\in\gamma(\A)$.
All three hexagons on Fig.~7 belong to one of the orbits. Hexagons of the other
orbit can be obtained from them by the $\pi/2$ rotation around the origin (the
eigenvectors of~$A$, which direct $l_1$ and $l_2$, are orthogonal, because $A^
{\T}=A$). As $k\to-\infty$, the hexagon $A^kW$ looks more and more like the
line~$l_1$; as $k\to\infty$, it looks more and more like~$l_2$. Directions of
these lines are points of the projective line $\R P^1=S^1$, which is the
absolute in the Poincar\'e circle model of the Lobachevskii plane. We encourage
the reader to find the relation between $\A$, $\gamma(\A)$, and the geodesic
line that connects these two absolute points. Hint: see~\cite{18,~\S17}
and~\S4.2 below.						\endexample

\subhead 2.3. Spines of torus bundles \endsubhead

From now on, $M$ denotes the total space of an orientable $T^2$-bundle over the
circle and $\A\in\SL(2,Z)$ is the monodromy operator (acting on the
one-dimensional homology group of the fiber containing the base point of~$M$)
of the bundle. By $M(\A)$ denote the manifold $M$ corresponding to the
monodromy operator~$\A$.

In this section, we construct a pseudominimal (see Definition~12 below) special
spine of $M(\A)$ with $\max(6,c(\A)+5)$ vertices. First, let us note that the
number $c(\A)$ is well defined by the manifold $M^3$ because of the following
statement.

\proclaim{Theorem~9~\cite{27}} Suppose that \rom3-manifolds $M(A)$ and $M(B)$,
where $A,B\in\SL(2,\Z)$, are homeomorphic. Then $B$ is $\GL(2,\Z)$-conjugate to
either $A$ or~$A^{-1}$. 					\endproclaim

If $A$ and $B$ are $\SL(2,\Z)$-conjugate, they yield the same number~$c(\A)$.
If $B\sim A^{-1}$, then $c(\Cal B)=c(\A)$ by virtue of~eq.~(3). Since the group
$\GL(2,\Z)$ is generated by its subgroup $\SL(2,\Z)$ and the element $C=
\pmatrix0&1\\ 1&0\endpmatrix$, it now suffices to check that conjugation by~$C$
does not affect the complexity~$c(A)$; however, this is clear from
Definition~10.

\definition{Definition~12~\cite{14}} A 2-dimensional component $\a$ of a
special polyhedron has a {\it counterpass\/} if its boundary $\dd\a$ passes
along some edge of $SP$ in both directions; it is called a {\it component with
short boundary\/} if $\dd\a$ passes through at most 3 vertices and visits any
of them only once. A special spine of a 3-dimensional manifold is said to be
{\it pseudominimal\/} if it contains neither components with counterpasses nor
components with short boundaries.				\enddefinition

\midinsert

\epsfxsize=360pt

\centerline{\epsffile{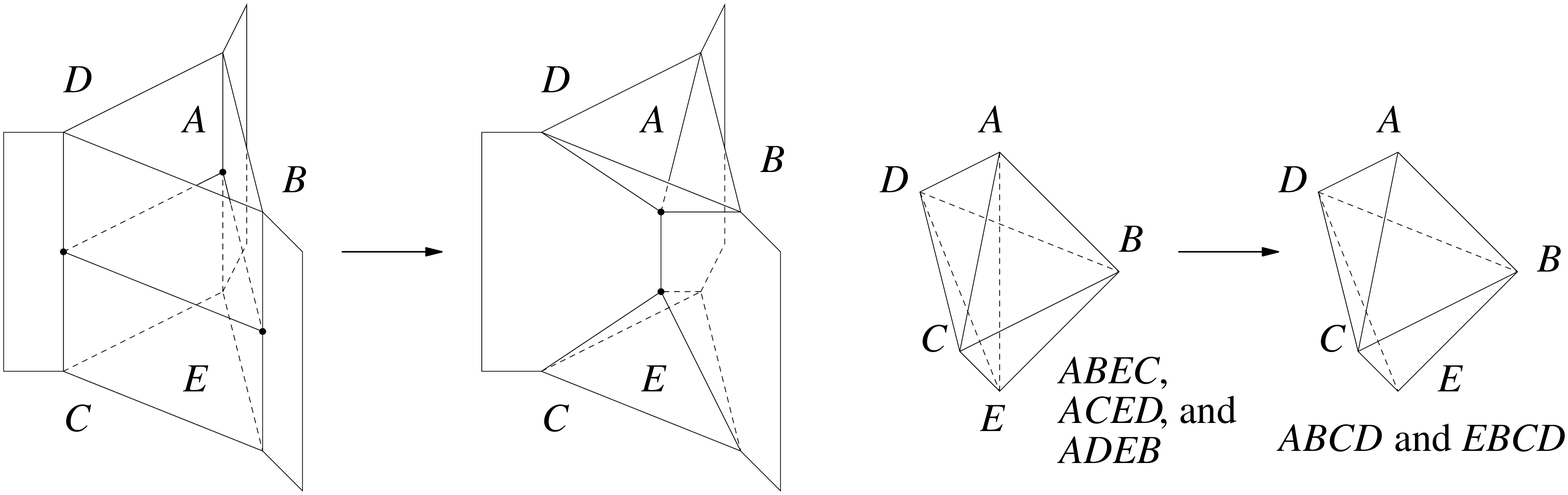}}

\botcaption{Figure 8} A simplification move (left) and the corresponding
Pachner move (right)	\endcaption

\endinsert

If a special spine $P$ is not pseudominimal, it is not minimal, because one can
apply simplification moves (see~\cite{14}) to~$P$ and get an almost simple
spine with a smaller number of vertices. For example, Figure~8 shows the effect
of a simplification move applied to a special spine with a triangular component
(the middle horizontal triangle in the left part of Fig.~8); it is easy to see
that the neighborhood of a 2-cell with short boundary of length~3 in a special
polyhedron~$P$ looks like the left hand side of Fig.~8. This move does not
change the spine outside of the fragment shown on Fig.~8. Note that the spine
obtained is special again: the move produces neither closed triple lines nor
non-cellular 2-dimensional components.

\remark{Remark} Consider the singular triangulation dual to a special spine
with a triangular component. Then the simplification move shown on Fig.~8
corresponds to the three-dimensional $(3,2)$ Pachner move~\cite{20}, which
replaces three tetrahedra by two tetrahedra. In the two-dimensional case, a
flip (see Fig.~3) corresponds to the $(2,2)$ Pachner move, which switches the
diagonal in a quadrilateral formed by two neighboring triangles. Recall that
the move shown on Fig.~8 and its inverse are sufficient to convert any two
special spines of the same compact three-dimensional manifold to one another,
see~\cite{12}; this fact is crucial for the construction of the Turaev--Viro
invariants~\cite{28}.						\endremark

\midinsert

\epsfxsize=185pt

\centerline{\epsffile{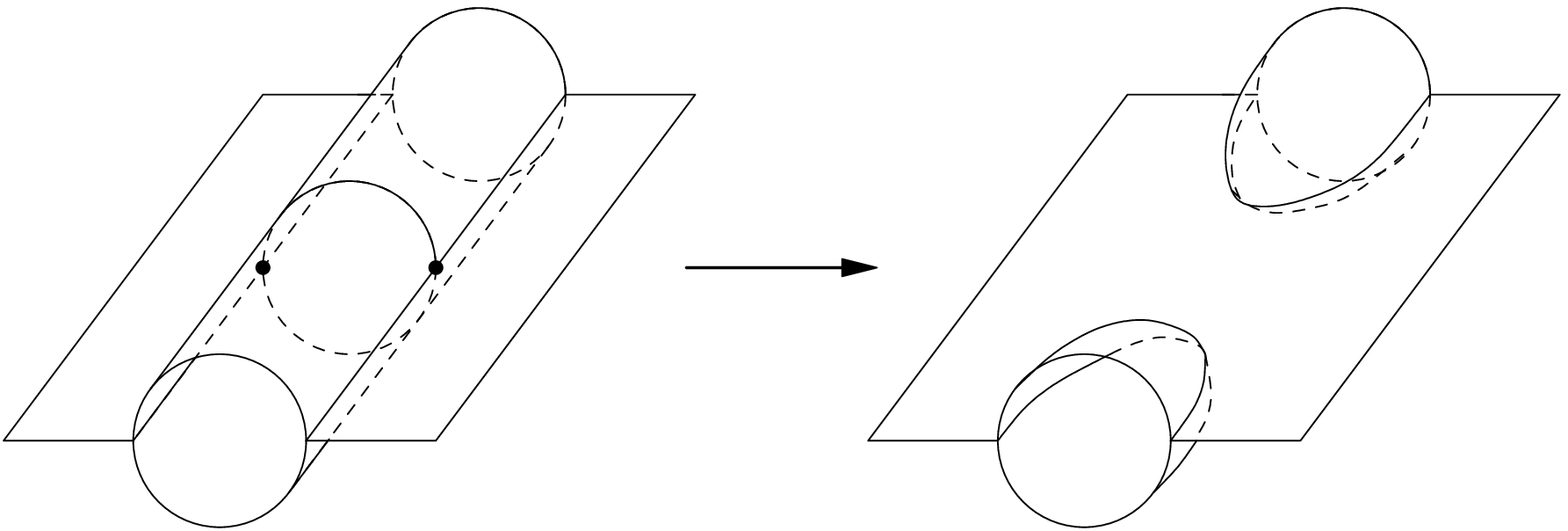} \epsffile{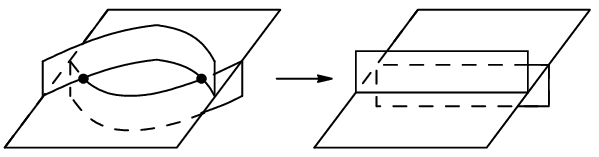}}

\botcaption{Figure 9} Two presentations of another simplification move
\endcaption

\endinsert

Figure~9 represents another simplification move, which is applicable to spines
containing a component with short boundary of length~2; clearly, the
neighborhood of this component looks like the left hand side of~Fig.~9. This
simplification move yields a simple, but not necessarily special, spine of the
same manifold (provided that the move had been applied to a simple spine).

To construct a spine of $M(\A)$, consider a fiber $T^2\x\{0\}$ and choose a
\thc $L_0$ in it; by doing so, we also fix a \thc $L_1$ in $T^2\x\{1\}$; note
that $W(L_1)=\A W(L_0)$. This choice is equivalent to the choice of some basis
in the lattice $H_1(T^2,\Z)$; by $A$ denote the matrix of $\A$ in this basis.
Construct a continuous family $L_t$ transforming $L_0$ into $L_1$ by isotopy
and $c(A)$ flips. Set $P_0=\bigcup\limits_{t\in[0,1]}L_t$; we assume that each
$L_t$ is embedded in $T^2\x\{t\}$. Note that $P_0$ is a simple polyhedron,
which is a spine of some punctured torus bundle. Two-dimensional components of
$P_0$ come from edges of $L_t$ as $t$ varies; similarly, one-dimensional
components of $P_0$ come from vertices of~$L_t$. The $c(A)$ flips correspond to
the vertices of~$P_0$, see~Fig.~10.

\midinsert

\epsfxsize=200pt

\centerline{\epsffile{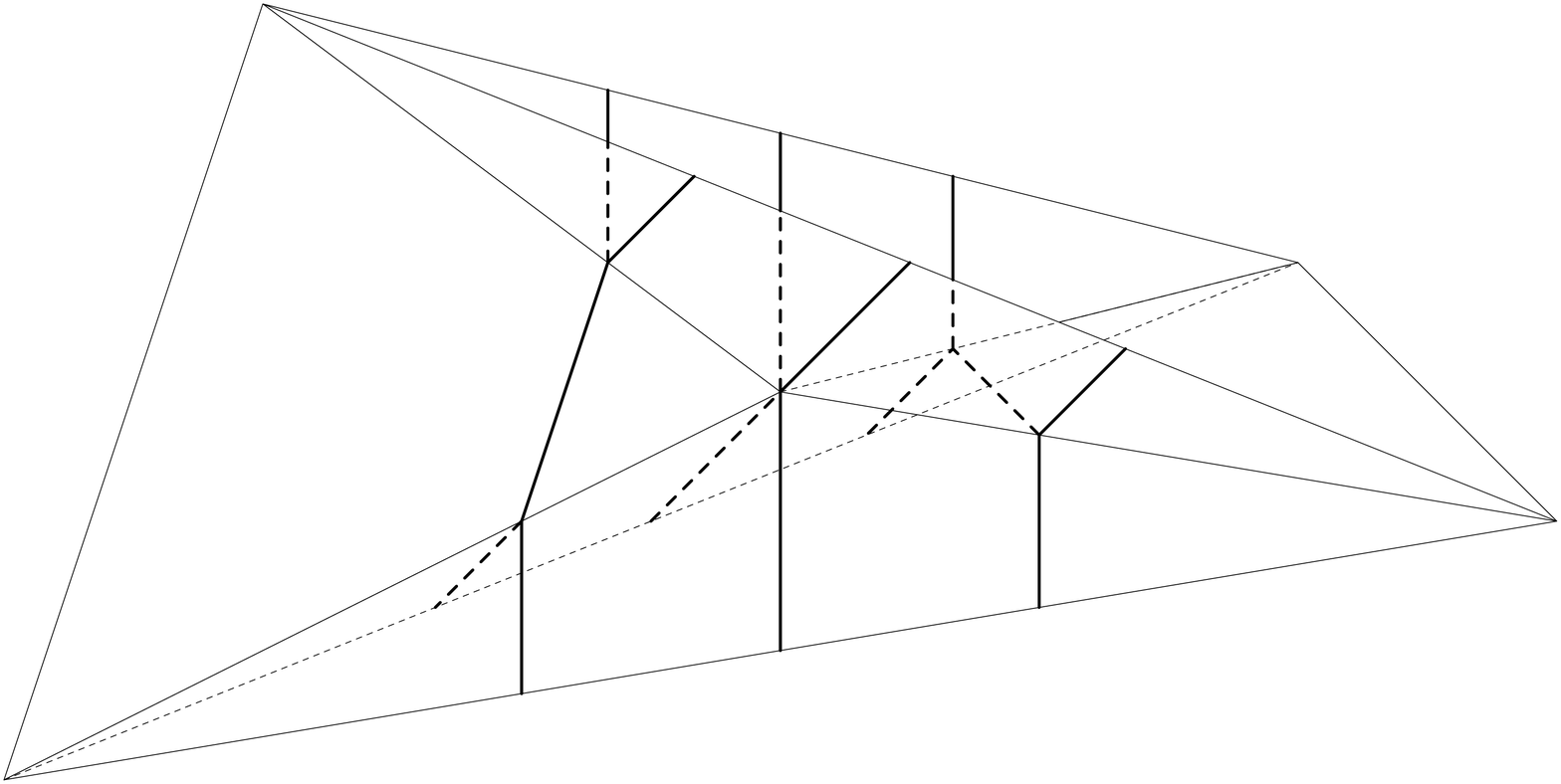}}

\botcaption{Figure 10} Vertices correspond to flips	\endcaption

\endinsert

To minimize the number of vertices, it is natural to choose a basis in $H_1(T^2
,\Z)$ so that the operator $A$ is minimal for~$\A$. If the operator $A$ is not
minimal, Theorem~7 guarantees that the last flip of the first round along the
base circle of the fibering and the first flip of the second round are mutually
inverse. This means that the second simplification move (see Fig.~9) is
applicable. Apply it until it is no longer possible. This process is nothing
but the construction of a minimal hexagon for $\A$ by the algorithm described
below the proof of Theorem~7. In the following, we suppose that $A$ is a
minimal matrix.

\example{Examples.~1} The only periodic operators $\A$ with $c(\A)>0$ have the
minimal matrices $\pmatrix\ph-0&\mp1\\ \pm1&\ph-0\endpmatrix$. This is a very
interesting case. The polyhedron $P_0$ constructed above has one vertex.
Consider the two-sheeted covering of the base $S^1$ of the fibering. It induces
the two-sheeted covering of the total space by the manifold $M(\A^2)=M(-I)$.
The preimage of $P_0$ under the covering is a polyhedron in $M(-I)$ with two
vertices, which can be cancelled by the second simplification move in two
different ways. This is the only (up to a sign and conjugacy) operator such
that $c(\A^k)<|k|c(\A)$ for $k\in\Z$, $|k|\ge2$, see Theorem~8; the other
periodic operators are of complexity~0.

{\bf2.} If $c(A)=0$, there are no flips at all. In this case $P_0$ contains no
vertices and consists of three orientable annuli and three edges if $A=I$, of
three nonorientable annuli and one edge if $A=-I$, of one nonorientable annulus
and one edge if $A$ is equal to the standard hexagon rotation matrix $R_{\pi/3}
=\pmatrix0&-1\\1&\ph-1\endpmatrix$ or its inverse (note that $R^6_{\pi/3}=I$),
and of one orientable annulus and two edges if $A$ equals $R^2_{\pi/3}$ or $R^4_
{\pi/3}$. This exhausts the case~$c(A)=0$.			\endexample

The polyhedron $P_0$ is not a spine of $M(\A)$, because the fibered space $M(\A
)$ admits a section that does not intersect $P_0$. This section represents a
nontrivial element of the group $\pi_1(M(\A))$, while the complement to a spine
of a closed manifold is a cell and hence cannot contain nontrivial loops. Let
us put $P_1=P_0\cup(T^2\x\{0\})$.

\proclaim{Lemma~2} $P_1$ is a spine of~$M(\A)$. 		\endproclaim

\demo{Proof} It is sufficient to show that $M(\A)\sm P_1$ is a 3-dimensional
cell. We have $M(\A)\sm P_1=T^2\x(0,1)\sm P_0=T^2\x(0,1)\sm\bigcup\limits_{t\in
(0,1)}L_t=\bigcup\limits_{t\in(0,1)}(T^2\x\{t\}\sm L_t)$, and Lemma follows.
								\qed\enddemo

Note that $P_1$ is not a simple polyhedron. Indeed, its part $T^2\x\{0\}$
contains a singular subset $L_0$, which is more complicated than a triple line:
three edges of $L_0$ yield three lines of transversal intersection of two
surfaces, and any of two vertices of~$L_0$ gives rise to a transversal
intersection of a triple line with one extra surface.

Let us modify the previous construction by gluing $T^2\x\{0\}$ with $T^2\x\{1
\}$ along a homeomorphism $\A+\delta$, where $\delta$ is a small shift of the
torus in a direction transversal to the edges of~$L_0$, see~Fig.~11. Put $P_2=
\bigcup\limits_{t\in[0,1]}L_t\cup(T^2\x\{0\})$. Again, $P_2$ is a spine
of~$M(\A)$.

\midinsert

\epsfxsize=100pt

\centerline{\epsffile{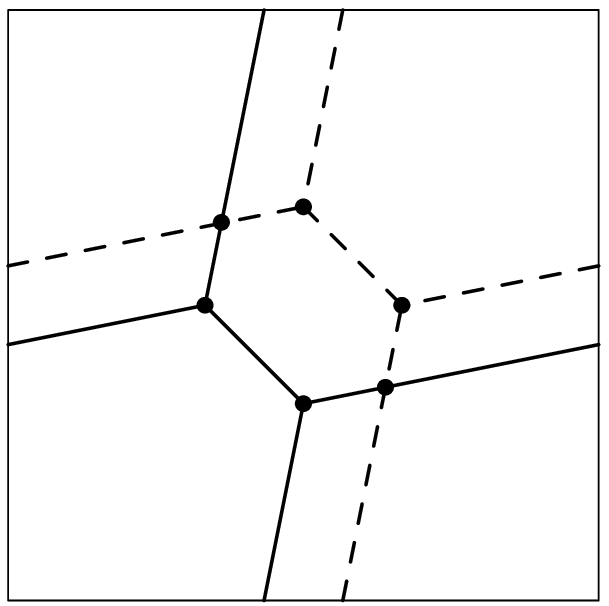}}

\botcaption{Figure 11} \thc $L_0$ and its $\delta$-shift $L_1$ (dashed line) in
the fiber $T^2\x\{0\}$ of~$M(\A)$				\endcaption

\endinsert

\proclaim{Lemma~3} $P_2$ is a special spine of~$M(\A)$ with $c(\A)+6$ vertices.
								\endproclaim

\demo{Proof} It is clear from the construction that $P_2$ is a simple
polyhedron. Its triple lines are the ``trajectories'' (as $t$ varies) of the
vertices of~$L_t$ and the ten segments of $L_0$ and $L_1$ shown on Fig.~11,
where the torus is represented by a square with the opposite sides to be
identified. There are $c(\A)$ vertices of $P_2$ that correspond to $c(\A)$
flips between $L_0$ at $t=0$ and $L_1$ at $t=1$, and six other vertices thar
are drawn on Fig.~11. Two of them arise from $T^2\x\{0\}$ and $L_t$, $0\le t<
\e$, and their neighborhoods in~$P_2$ look like Fig.~1\,d. Two other vertices
on Fig.~11 arise from $T^2\x\{0\}=T^2\x\{1\}$ and $L_t$, $1-\e<t\le1$; their
neighborhoods look like the horizontal mirror reflection of Fig.~1\,d. The last
two vertices on Fig.~11 correspond to two intersection points of $L_0$ and $L_
1$, and their neighborhoods look like Fig.~1\,c. Thus $P_2$ is a simple spine
of $M(\A)$ with $c(\A)+6$ vertices.

It remains to prove that $SP_2$ contains no closed triple lines and all
connected components of $P_2\sm SP_2$ are 2-dimensional cells, cf\.
Definition~2. First group of triple lines of $SP_2$ is formed by ten arcs in
$T^2\x\{0\}$ shown on Fig.~11. Obviously, they are not closed. The rest $2c(A)+
2$ triple lines are swept by the vertices of the \thcs $L_t\subset T^2\x\{t\}$,
$0<t<1$. They end at vertices of~$P_2$, too, and thus are not closed.

Connected components of $P_2\sm SP_2$ also belong to two groups. Four of them,
two hexagonal and two quadrilateral, lie in the fiber $T^2\x\{0\}$, see Fig.~11.
They are cells. Any other connected component of $P_2\sm SP_2$ intersects any
fiber $T^2\x\{t\}$, $a<t<b$ (where $a$ is equal to either $0$ or one of the
flip moments, and $b$ is either one of the flip moments or is equal to~$1$),
along one edge of $L_t$, and does not intersect other fibers; this implies that
this component is a cell. We have proved that the polyhedron $P_2$ is special.
								\qed\enddemo

\proclaim{Corollary} $c(M(\A))\le c(\A)+6$.		\qed\endproclaim

\example{Example} Three-dimensional torus can be represented as~$M(I)$, $I=
\pmatrix 1&0\\ 0&1\endpmatrix$. Since $c(I)=0$, the construction above gives a
special spine of~$T^3$ with six vertices. The manifold $T^3=M(I)$ is contained
in Table~7 of the preprint~\cite{14} under the name~$6_{71}$. It is shown
in~\cite{14} that all manifolds of complexity at most~5 are different from~$T^
3$. So we have $c(T^3)=6$. The spine that we constructed here does not differ
from the spine $6_{71}$ from~\cite{14,~\S 5.2}, while our way of presenting
spines differs significantly from the one used in~\cite{14}.

\midinsert

\epsfysize=100pt

\centerline{\epsffile{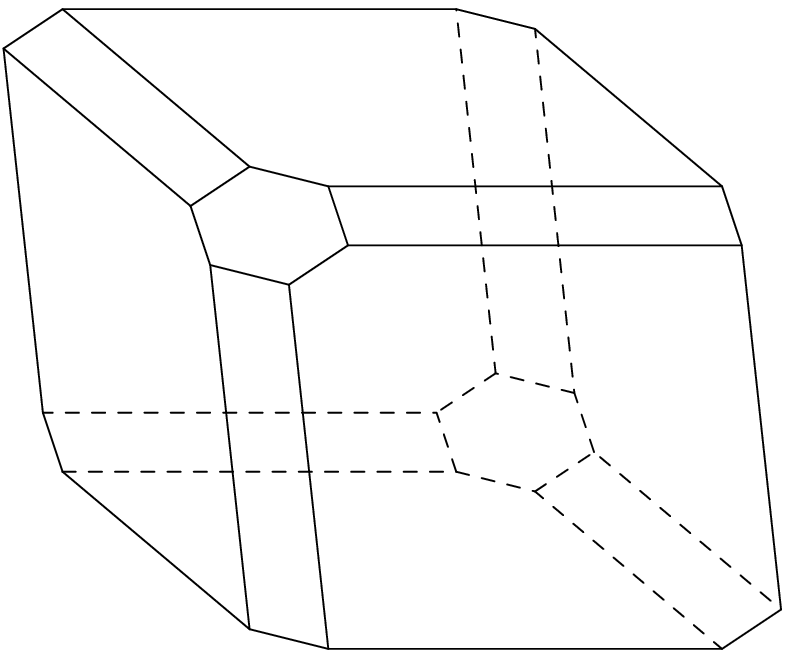}}

\botcaption{Figure 12} The complement $T^3\sm P_2$ of the minimal spine~$P_1$
of~$T^3$	\endcaption

\endinsert

The torus $T^3$ can be obtained from the cube by gluing its opposite faces.
This yields a natural cell decomposition of~$T^3$ with one vertex, three edges,
three ``square'' 2-dimensional cells and one 3-dimensional cell. The
2-dimensional skeleton $\sk_2(T^3)$ has singular points more complicated than
triple lines and vertices of simple polyhedra. However, the minimal spine
of~$T^3$ can be obtained as a small perturbation of~$\sk_2(T^3)$.

Let the \thcs $L_t$, $t\in[0,1]$, be very close to the bouquet of a parallel
and a meridian of~$T^2\x\{t\}$, and let the shift~$\delta$ involved in the
construction of~$P_2$ be very small. Then the 3-dimensional cell $T^3\sm P_2$
is very close to the 3-dimensional cube. Figure~12 represents this cell. If we
identify opposite faces of this polyhedron by parallel transports (or,
equivalently, tessellate~$\R^3$ into parallel copies of this polyhedron and
consider a quotient over the appropriate lattice~$\Z^3$), we get the torus~$T^
3$; the image of the boundary of the polyhedron under this gluing is the
minimal spine of~$T^3$ close to~$\sk_2(T^3)$.			\endexample

The same construction gives special spines with six vertices for the manifolds
$M(-I)=6_{70}$,
	$$M\left(\pmatrix -1&-1\\ \ph-1&\ph-0\endpmatrix\right)=
		M\left(\pmatrix\ph-0&\ph-1\\ -1&-1\endpmatrix\right)=6_{67},$$ 
and 
	$$M\left(\pmatrix 0&-1\\ 1&\ph-1\endpmatrix\right)=M\left
			(\pmatrix \ph-1&1\\ -1&0\endpmatrix\right)=6_{65}.$$
The spines constructed in this way are minimal spines of these manifolds,
because all of them are of complexity~6; in fact, all manifolds of complexity
up to~5 are quotient spaces of the sphere~$S^3$, see~\cite{14}.

However, in all other cases (that is, if $c(\A)>0$) the spines with $c(\A)+6$
vertices are not minimal spines of the manifolds~$M(\A)$. For example, the
spaces
$$6_{66}=M\left(\pmatrix\ph-0&1\\ -1&0\endpmatrix\right),\quad
  6_{68}=M\left(\pmatrix-1&\ph-0\\ -1&-1\endpmatrix\right),\quad\hbox{and}\quad
  6_{69}=M\left(\pmatrix 1&0\\ 1&1\endpmatrix\right)$$
are manifolds of complexity~6, while
$$c\left(\pmatrix\ph-0&1\\ -1&0\endpmatrix\right)=
  c\left(\pmatrix -1&\ph-0\\ -1&-1\endpmatrix\right)=
  c\left(\pmatrix 1&0\\ 1&1\endpmatrix\right)=1,$$
and our construction gives their spines with 7 vertices.

This happens because some of the spines with $c(\A)+6$ vertices constructed
above are not pseudominimal whenever $c(\A)>0$. Namely, they have a triangular
component, and the first simplification move (see Fig.~8) can be applied.

Let us return to Fig.~11. Assume that the first flip in the sequence taking
$L_0$ to $L_1$ involves the short edge of~$L_0$, that is, the edge that does
not intersect dashed lines on Fig.~11. This condition can be satisfied by an
appropriate choice of the shift~$\delta$ involved in the construction of~$P_2$.
Then the 2-dimensional cell of $P_2$ adjacent to this edge and not contained in
$T^2\x\{0\}$ is a triangle, and we can apply the first simplification move,
which gives a spine of $M(\A)$ with a smaller number of vertices. This spine
can be described in other words as follows. Let $L'$ be the \thc obtained after
the first of $c(\A)$ flips converting $L_0$ into~$L_1$. Glue the square from
Fig.~13 into the torus $T^2\x\{0\}\in M(\A)$ and embed $L_t$ in $T^2\x\{t\}$
for all $t\in(0,1)$, where the family $L_t$ contains $c(\A)-1$ flips and
connects $L'$ with~$L_1$. Note that the first of $c(\A)-1$ flips converting
$L'$ into $L_1$ is performed along a long edge of $L'$, because a flip along
the short edge would annihilate with the flip converting $L_0$ to~$L'$. The
new spine $P_3$ has $6+c(\A)-1=c(\A)+5$ vertices. So we have $c(M(\A))\le c(\A)
+5$ whenever $c(\A)>0$.

\midinsert

\epsfxsize=100pt

\centerline{\epsffile{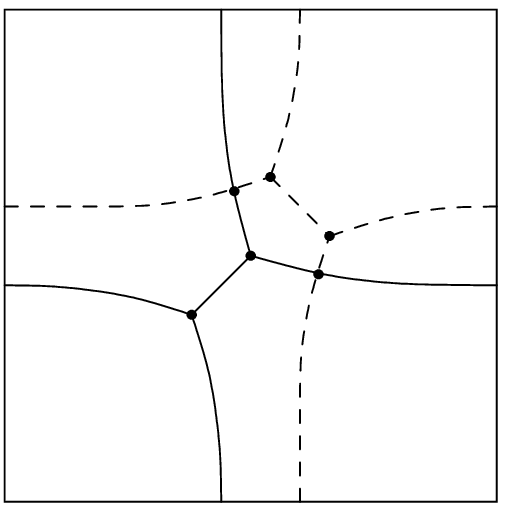}}

\botcaption{Figure 13} \thcs $L'$ and $L_1$ (dashed) in $T^2\x\{0\}$
\endcaption

\endinsert

The proof that the spine $P_3$ is special repeats the proof of Lemma~3.

\proclaim{Theorem~10}

\rom{1)} $c(M(\A))\le\max(6,c(\A)+5)$.

\rom{2)} The spine $P_3$ constructed above is pseudominimal.	\endproclaim

Compare the first statement with Corollary of Lemma~3.

\demo{Proof} If $c(\A)=0$, then $c(M(\A))=6$. So we may suppose that $c(\A)>0$.
Since $P_3$ is an almost simple (and even special) spine of $M(\A)$ with $c(\A)
+5$ vertices, the first statement is obvious. The argument similar to the proof
of Lemma~3 shows that two-dimensional cells of $P_3$ have no counterpasses.
So we only have to show that $P_3$ has no components with short boundaries. The
four cells contained in $T^2\x\{0\}$ are pentagons, see Fig.~13. The cells that
have no boundary edges in $T^2\x\{0\}$ have even numbers of edges, namely, $2k-
2$, where $k$ is the number of flips from the vertex where the cell appears to
the vertex where the cell disappears (including both the first flip and the
last one). The matrix $A$ involved in the construction of $P_3$ is a minimal
matrix of~$\A$. This implies that any two consecutive flips in the sequence
involved in the construction are not inverse to one another, that is, $k>2$ for
any cell considered above.

It remains to consider at most 6 two-dimensional cells that have an edge in $T^
2\x\{0\}=T^2\x\{1\}$ (if $A$ is, up to a sign, a power of a Jordan block, there
are only 5 cells of this type; otherwise, no 2-cell touches both $T^2\x\{0\}$
along an edge of $L'$ and $T^2\x\{1\}$ along an edge of $L_1$, so three edges
of $L'$ and three edges of $L_1$ belong to six different cells of $P_3\sm T^2\x
\{0\}$). Two cells of $P_3\sm T^2\x\{0\}$ are adjacent to the long edges of $L'$
(the edges that intersect dashed lines on Fig.~13). Each of these cells has at
least 4 boundary edges: two segments of a long edge of $L'$ and two edges
(transversal to fibers) that arise from the vertices of~$L'$. The same argument
works for two cells of $P_3\sm T^2\x\{0\}$ adjacent to the long edges of~$L_1$.
Consider the cell of $P_3\sm T^2\x\{0\}$ adjacent to the short edge of~$L'$. Of
course, it has at least 3 edges: the short edge of $L'$ and two edges that are
trajectories of the vertices of $L'$ as $t$ varies. It has at least one more
edge: otherwise, the first flip in the sequence of flips converting $L'$ to $L_
1$ is performed along the short edge of $L'$, which is impossible by the
construction of $P_3$, see above. For the same reason, the last flip (which
results in $L_1$) cannot be performed along the short edge of~$L_1$: by the
minimality of the matrix $A$, it cannot be cancelled with the flip connecting
$L_0$ with $L'$, see Theorem~7. This means that the 2-dimensional cell of $P_3
\sm T^2\x\{1\}$ adjacent to the short edge of $L_1$ also has more than 3 edges,
and $P_3$ contains no components with short boundaries. The Theorem is proved.
								\qed\enddemo

\proclaim{Conjecture~2} The pseudominimal spines of the manifolds $M(\A)$
constructed above are in fact their minimal spines, and the upper bound for
complexity given in Theorem~\rom{10} is in fact its exact value\rom: $c(M(\A))=
\max(6,c(\A)+5)$ for any monodromy operator $\A\in\SL(2,\Z)$. In other words,
any singular triangulation of $M(\A)$ involves at least $c(\A)+5$ tetrahedra if
$c(\A)>0$ and \rom6 tetrahedra if $c(\A)=0$.			\endproclaim

\subhead 2.4. Digression: spines of lens spaces 	\endsubhead

Pseudominimal special spines of the lens spaces~$L_{p,q}$, $p>3$, with exactly
$E(p,q)-3$ vertices were constructed in~\cite{14}. In that paper, spines are
presented by drawing the neighborhood of the singular graph of a spine. This
allows to draw spines on the plane; however, it remains unclear how the spines
are embedded into corresponding manifolds.

In this section, we construct pseudominimal special spines of~$L_{p,q}$, $p>3$,
with $E(p,q)-3$ vertices, making use of the techniques developed in~\S2. We
omit some details and proofs.

Consider two solid tori. The meridians of their boundary tori are well defined,
while the parallels are defined modulo meridians only. Let $\mu_0,\mu_1$ be the
meridians of the tori and $\sigma_0,\sigma_1$ be their parallels such that the
pair of the oriented cycles $(\sigma_0,\mu_0)$ defines the positive orientation
of the boundary of the first torus and the pair $(\sigma_1,\mu_1)$ defines the
negative orientation of the boundary of the second torus. There is a unique
pair of positive integer numbers $(r,s)$ such that $r<p$, $s<p$, and $qs-pr=1$.
Put $A=\pmatrix s&p\\ r&q\endpmatrix$ and attach the solid tori to one another
so that the induced homomorphism of the one-di\-men\-si\-on\-al homology
groups of their boundary tori has the matrix~$A$ (in the bases $(\sigma_0,\mu_
0)$ and $(\sigma_1,\mu_1)$). We get a closed orientable 3-manifold that is
nothing but~$L_{p,q}$.

Note that $A\in\SL(2,\Z)$, $c(A)=E(p,q)$, and the parallels $\sigma_0$ and
$\sigma_1$ represent nontrivial elements of $\pi_1(L_{p,q})=\Z_p$. This implies
that any spine of $L_{p,q}$ intersects these loops. Let us shift $\sigma_0$
in the interior of the first solid torus and consider the tubular neighborhood
$U_0$ of the shifted curve. Obviously, $U_0$ is a solid torus. Similarly,
construct $U_1$ as a tubular neighborhood of $\sigma_1$ shifted inside of the
interior of the second torus. We may assume $U_0$ and $U_1$ to be disjoint.
Then $L_{p,q}=U_0\cup (T^2\x[0,1])\cup U_1$. Let $L_i$, $i=0,1$, be standard
(with respect to the bases $(\sigma_i,\mu_i)$) \thcs in the tori $T^2_i=\dd U_
i$; they are defined up to isotopy. Following the construction of \S2.3,
consider a continuous family $L_t\subset T^2\x\{t\}$ connecting $L_0$ to $L_1$
with $c(A)$ flips. Put $P_0=\bigcup\limits_{t\in[0,1]}L_t$. Let $D_i$, $i=0,1$,
be meridional disks of the $U_i$ intersecting $L_i$ transversally at one point.
Put $P_1=D_0\cup T^2_0\cup P_0\cup T^2_1\cup D_1$.

\proclaim{Lemma~4} The polyhedron $P_1$ is a special spine of $L_{p,q}$ with
three punctures. It has $E(p,q)+6$ vertices.			\endproclaim

\demo{Proof} The complement $L_{p,q}\sm P_1$ consists of three cells $U_0\sm D_
0$, $(T^2\x[0,1])\sm P_0$, and $U_1\sm D_1$. There are $c(A)=E(p,q)$ vertices
in the interior part of $P_0$. Further, there are 3 vertices on $T^2_0$, which
correspond to two vertices of $L_0$ and the intersection point of $L_0$
and~$\dd D_0$. Similarly, there are 3 vertices of $P_1$ on~$T^2_1$. It remains
to show that $P_1$ is a special polyhedron. This can be proven by analogy with
Lemma~3.							\qed\enddemo

Below we show that one can decrease the number of vertices ``inside of $P_0$\,''
by one and the number of vertices ``near each $U_i$'' by four. This gives a
spine with $E(p,q)+6-1-4-4=E(p,q)-3$ vertices.

Recall that the parallels $\sigma_i$ are defined only modulo meridians~$\mu_i$.
Thus, the \thcs $L_i$ are defined only up to powers of the Dehn twists along
the meridians, that is, up to transformations $\sigma_i\mapsto\sigma_i+n_i\mu_
i$, $n_i\in\Z$. By varying $n_0$ and $n_1$, one can decrease the distance
in~$\Gamma$ between $B^{n_0}W_0$ and $C^{n_1}\!AW_0$ and thus decrease the
number of the vertices inside of~$P_0$; here $W_0$ is the standard hexagon, $B=
\pmatrix1&0\\1&1\endpmatrix$ is the matrix representing the Dehn twist
along~$\mu_0$, and $C=ABA^{-1}$ is the matrix corresponding to the Dehn twist
along~$\mu_1$.

\proclaim{Lemma~5} $\min\limits_{n_0,n_1\in\Z}d(B^{n_0}W_0,C^{n_1}\!AW_0)=E(p,q
)-1$.								\endproclaim

\demo{Proof} Note that the graph $\Gamma$ contains edges $B^uW_0\,B^{u+1}W_0$
and $C^v\!\!AW_0\,C^{v+1}\!AW_0$ for all $u,v\in\Z$. Since $\Gamma$ is a tree,
there are $m_0,m_1\in\Z$ such that the path from $B^{n_0}W_0$ to $C^{n_1}\!AW_
0$ for all $n_0,n_1\in\Z$ consists of the following three legs: $B^{n_0}W_0\,B^
{m_0}W_0$, $B^{m_0}W_0\,C^{m_1}\!AW_0$, and $C^{m_1}\!AW_0\,C^{n_1}\!AW_0$. Now
it is obvious that $\min\limits_{n_0,n_1\in\Z}d(B^{n_0}W_0,C^{n_1}\!AW_0)=d(B^{
m_0}W_0,C^{m_1}\!AW_0)$. By considering the three legs of the path from $W_0$
to $AW_0$, one can see that $m_0=1$ (because $p>q>0$), $m_1=0$ (because both
positive and negative Dehn twists along $\mu_1$ do not affect the leading
vertex $(p,q)$ of $AW_0$ and thus  increase the distance to $W_0$), and the
length of the middle leg of this path is $d(BW_0,AW_0)=d(W_0,AW_0)-|m_0|-|m_1|=
E(p,q)-1$. Also see Theorem~17 in~\S4.3.			\qed\enddemo

\midinsert

\epsfxsize=369pt

\centerline{\epsffile{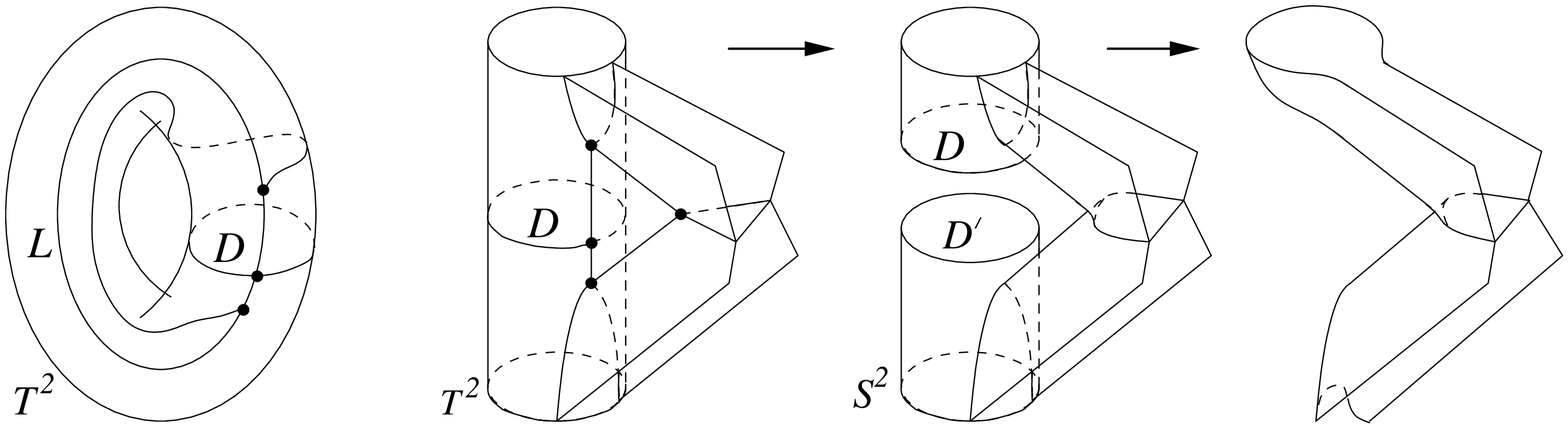}}

\botcaption{Figure 14} Simplification of $P_0$ near $T^2_i$
\endcaption

\endinsert

By virtue of Lemma~5, we can decrease by one the number of the vertices inside
of $P_0$ by another choice of a $\th$-curve~$L_0$. Now we have a special spine
of $L_{p,q}$ with three punctures having $E(p,q)+5$ vertices. The disks $D_0$
and $D_1$ are components with short boundaries. By $\tau_i$ denote the edge
of $L_i$ that intersects~$\dd D_i$. Note that two other edges of $L_i$ form the
meridian of $T_i$, the first flip in the sequence connecting $L_0$ with $L_1$
is performed along $\tau_0$ while the last flip in this sequence is performed
along~$\tau_1$ (flips along other edges are equivalent to meridional Dehn
twists and thus do not lead out of the mainstreams $\gamma(B)=\{B^{n_0}W_0\mid
n_0\in\Z\}$ and $\gamma(C)=\{C^{n_1}\!AW_0\mid n_1\in\Z\}$, while the path
between $L_0$ and $L_1$ is the shortest path that connects these mainstreams).
We can apply the following simplification move in the neighborhood of~$D_i$.
First, add a parallel copy $D_i'$ of $D_i$. Second, delete the lateral surface
of the cylinder bounded by $D_i$, $D_i'$, and a thin strip of~$T^2_i$. Finally,
delete the cell of $P_1$ adjacent to $\tau_i$; this cell is triangular, because
$\tau_i$ is the edge involved in the flip in $P_0$ closest to $T^2_i$, see
Fig.~14. So, the first step adds one vertex on each $T^2_i$, the second step
kills two vertices on each~$T^2_i$, and the last step kills three vertices near
each of~$T^2_i$. By $P_1$ denote the polyhedron obtained by the construction
above. Obviously, it has $E(p,q)-3$ vertices. Further, one can see that two
remaining edges of $L_i$ (which differ from~$\tau_i$) form a closed triple
line~$S^1_i$ and the complement $L_{p,q}\sm P_1$ still consists of three cells,
two of which are bounded by the spheres $S^2_i$ obtained from the torus $T^2_i$
and their meridional disks $D_i$, $D_i'$ by deleting the thin strip bounded by
$\dd D_i$ and $\dd D_i'$ from~$T^2_i$. The circles $S^1_i$ divide the spheres
$S^2_i$ into two disks each; one of the disks contains $D_i$, the other
contains~$D_i'$. Delete from~$P_1$ the disks of the $S^2_i$ that contain~$D_i
'$. This yields a polyhedron $P$ with $E(p,q)-3$ vertices such that $L_{p,q}\sm
P$ is a cell.

\proclaim{Theorem~11} The polyhedron $P$ is a pseudominimal special spine of $L
_{p,q}$ with $E(p,q)-3$ vertices. It coincides with the spine of $L_{p,q}$
presented in~\cite{14}.					\qed\endproclaim

It was shown in~\cite{14} that the spines constructed in that paper are
pseudominimal. So it is sufficient to prove only the second statement; we leave
it to the reader.

\head\S3. Lower bound for $C^1$-smooth spines transversal to fibers \endhead

In this section, we prove the following statement.

\proclaim{Theorem~12} Let $\A$ be a non-periodic operator. Then for any $C^
1$-smooth spine $P$ of $M(\A)$ transversal to the fibers, we have $c(P)\ge
\dfrac15\,c(\A)+2$, where $c(P)$ stands for the number of the vertices of~$P$.
								\endproclaim

Recall that $M(\A)$ is the total space of the fibering $p\:M^3\overset{T^2}\to
\tto S^1$ with monodromy operator~$\A$.

\subhead 3.1. Morse transformations in simple polyhedra \endsubhead

Projection $p\:M(\A)\to S^1$ defines an $S^1$-valued function on~$M(\A)$. Let
us explore its restriction $p|_P$ to a simple spine $P$ of~$M(\A)$. The
polyhedron $P$ is a stratified space; its strata are the vertices, the edges,
and the two-di\-men\-si\-onal components of~$P$. We may assume that $p|_{\ol
\sigma}$ has an everywhere continuous derivative for any edge or
2-com\-po\-nent~$\sigma$: a $C^0$-small deformation of the embedding mapping
$i\:P\emb M^3$ is sufficient. Since the space of Morse functions on a manifold
is $C^1$-dense in the set of all $C^1$-smooth functions on it~\cite{16}, we may
assume that $p|_{\ol\sigma}$ is a Morse function for any edge or 2-component
$\sigma$ of~$P$. In this case $p$ is called a {\it Morse function\/} on~$P$.

Consider a fiber $F_t=p^{-1}(t)$ of the fibering $p\:M(\A)\to S^1$, where $t$
is a local coordinate on $S^1$. If $t$ is a nonsingular value of~$p$ (that is,
$F_t$ contains neither vertices of~$P$ nor critical points of the restrictions
of~$p$ to the edges and 2-com\-po\-nents of~$P$), then $F_t$ is transversal to
all strata of~$P$ and the intersection $K_t=F_t\cap P$ is a trivalent graph,
possibly disconnected. We will explore how $K_t$ changes as $t$ varies. For
much more comprehensive exposition of Morse transformations in stratified
spaces, see~\cite9.

Without loss of generality, it can be assumed that any singular fiber (that is,
a singular level of~$p$) contains only one singular point of $p|_P$ (including
vertices), and there are only finite number of singular fibers. Also, we may
assume that the restrictions $p|_{\ol\sigma}$ to the closures of triple lines
and two-di\-men\-si\-onal components of~$P$ have no boundary critical points.
Let two nonsingular fibers $F_-$ and $F_+$ be the preimages of two close points
$t_-,t_+\in S^1$. If the interval $(t_-,t_+)$ contains no singular values
of~$p$, then the graphs $K_-$ and $K_+$ are isotopic, that is, there is a
continuous family of embeddings $i_t\:K\emb F_t$, $t\in[t_-,t_+]$, of the same
graph~$K$ to the tori~$F_t$ such that $K_-=i_{t_-}(K)$ and $K_+=i_{t_+}(K)$.

\midinsert

\epsfxsize=200pt

\centerline{\epsffile{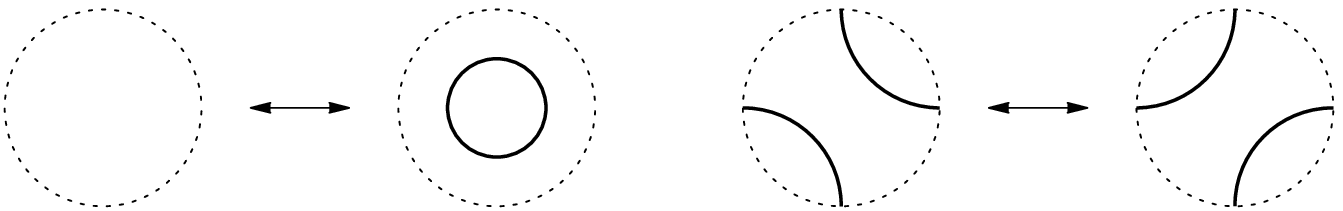}}

\botcaption{Figure 15} Minimax (left) and saddle (right) Morse transformations
\endcaption

\endinsert

Now suppose that there is exactly one singular fiber $F_0$ between $F_-$ and $F
_+$; by $t_0$ denote the corresponding critical value. Thus $F_0$ contains
either a singular point of the restriction of $p$ to a two-di\-men\-si\-onal
component of~$P$ or a singular point of the restriction of $p$ to a triple line
of~$P$, or a vertex of~$P$. In the first case, the difference between $K_-$ and
$K_+$ is nothing but the Morse transformation of the level set $Q=t-t_0$ of a
real quadratic form $Q(x,y)=\pm x^2\pm y^2$, see Fig.~15.

\midinsert

\epsfxsize=372pt

\centerline{\epsffile{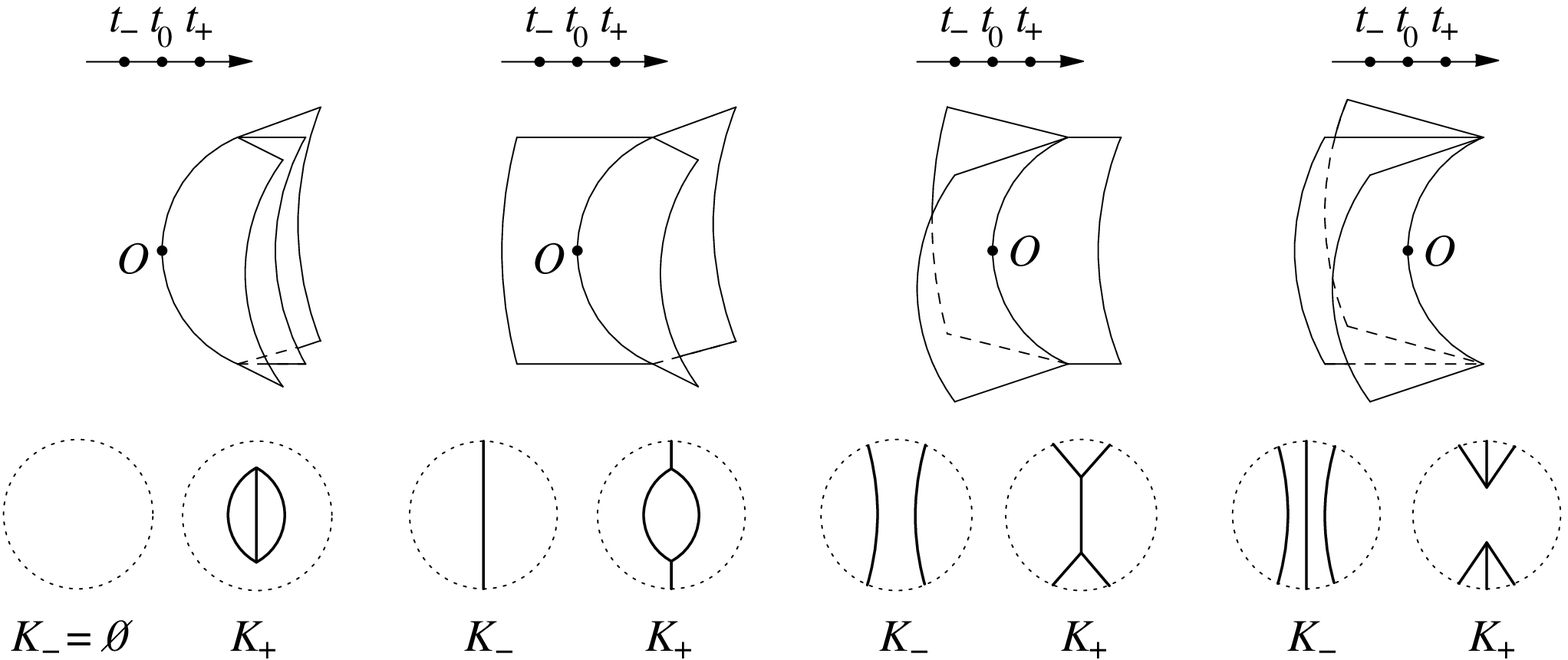}}

\botcaption{Figure 16} Transformation of $K$ induced by a minimum point on a
triple line							\endcaption

\endinsert

To explore the second case, suppose that the critical point of the restriction
of~$p$ to a triple line of $P$ is a minimum point of~$p$; denote this point
by~$O$. By $v_i$, $i=1,2,3$, denote the unit vector tangent to the $i$\snug th
two-di\-men\-si\-onal component $\sigma_i$ at~$O$ and orthogonal to the triple
line. Since $O$ is not a critical point of $p|_{\ol\sigma_i}$, we have $p_*v_i
\ne0$. Here the $\sigma_i$ are the components of $(P\sm SP)\cap U(O)$ of the
two-di\-men\-si\-onal strata of $P$ in the neighborhood $U(O)$ of~$O$; it does
not matter whether the $\sigma_i$ actually belong to different components of $P
\sm SP$. Let us say that the component $\sigma_i$ is {\it ascending\/} at~$O$
if the vector $p_*v_i\in T_{t_0}S^1$ induces a positive orientation of~$S^1$
and {\it descending\/} if this orientation is negative. There may be 0, 1, 2 or
3 descending components. Their number defines the transformation between $K_-$
and $K_+$ up to isotopy. The graphs $K_-$ and $K_+$ differ only inside of the
dotted circle, see Fig.~16. If $O$ is a maximum point of $p$ restricted to a
triple line, one has to revert all arrows and signs in the lower indices on
Fig.~16.

The last case (where the fiber $F_0$ contains a vertex $V\in P$) is more
complicated. We will return to it later.

\remark{Remark} The double of a disk is~$S^2$. Draw the fragment $K_-\cap U(O)$
of the graph $K_-$ (that is, the picture inside of the left dotted circle, see
Figs\. 15,~16) on the lower hemisphere of $S^2$ and the fragment $K_+\cap U(O)$
of the graph $K_+$ (that is, the picture inside of the right dotted circle) on
the upper hemisphere. In all four cases, we get an embedding of the link of~$O$
in~$P$ to the sphere with the equator drawn on it by a dotted line. The number
of intersection points of the equator and the link of~$O$ is twice the number
of descending components at~$O$.				\endremark

This is not a mere coincidence but a general recipe for describing the
difference between $K_-$ and~$K_+$. Suppose that there is exactly one singular
point~$O$ between $F_-$ and~$F_+$. Consider an $\e$-neigh\-bor\-hood $U(O)$
of~$O$, where $\e<t_0-t_-$ and $\e<t_+-t_0$. The intersection $P\cap U(O)$ is
the cone over the graph $\lk O$, which is the circle with none, two or three
radii, see the definition of a simple polyhedron (Definition~1 in~\S1). By $D_
-$, $D_0$, and $D_+$ denote the intersections of $U(O)$ with $F_-$, $F_0$, and
$F_ +$, respectively. Obviously, $K_+\sm D_+$ is isotopic to $K_-\sm D_-$.
Thus, the difference between $K_+$ and $K_-$ is ``hidden inside of~$U(O)$''.
The singular fiber $F_0$ cuts the sphere $S^2=\dd U(O)$ into two hemispheres
$S^2_-$ and $S^2_+$. Put $K_-'=(K_0\sm U(O))\cup (P\cap S^2_-)$ and $K_+'=(K_0
\sm U(O))\cup (P\cap S^2_+)$. These graphs are embedded into tori $(F_0\sm U(O)
)\cap S^2_-$ and $(F_0\sm U(O))\cap S^2_+$.

\proclaim{Lemma~6} The graph $K_-$, respectively, $K_+$, is isotopic to the
graph $K_-'$, respectively, $K_+'$.				\endproclaim

\demo{Proof} Let us define a surface $G_s$ as the union of $F_{t_0+s(t_0-t_-)}
\sm U(O)$ and $\{x\in S^2_-\mid p(x)<t_0+s(t_0-t_-)\}$, where $-1\le s\le0$,
and put $K_s'=P\cap G_s$. In $U(O)\sm O$, the edges and 2-com\-po\-nents of~$P$
are transversal to all fibers and to the spheres centered at~$O$. So all $G_s$,
$-1\le s\le0$, are transversal to~$P$, and the family $K_s'$ provides an
isotopy between $K'_{-1}\subset G_{-1}$ and $K'_0\subset G_0$. It remains to
note that $G_{-1}=F_-$, $K'_{-1}=K_-$, and $K'_0=K'_-$. An isotopy between $K_
+$ and $K_+'$ can be constructed in a similar way.		\qed\enddemo

\midinsert

\epsfxsize=372pt

\centerline{\epsffile{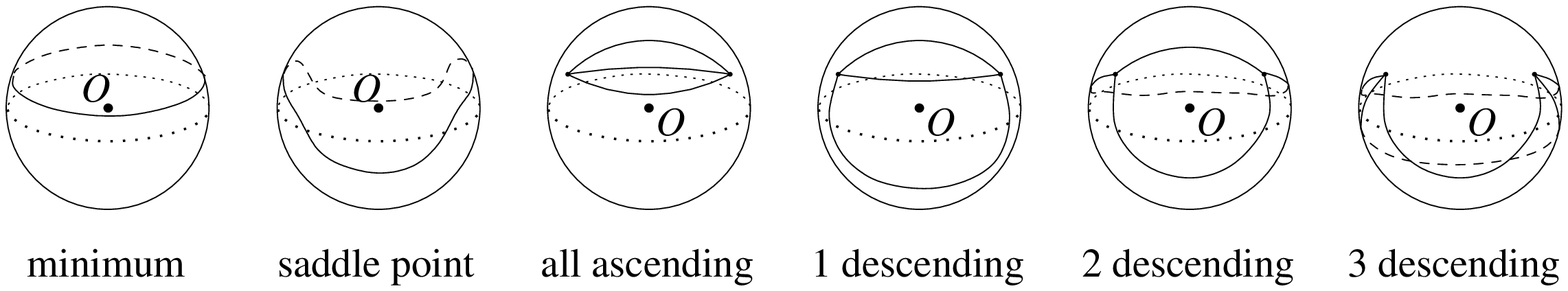}}

\botcaption{Figure 17} Mutual position of $\lk_{F_0}(O)$ (dotted line) and
$\lk_P(O)$ in $S^2$						\endcaption

\endinsert

So to describe a ``surgery'' in the graph $K$ near a singular point $O$, it is
sufficient to draw the links of~$O$ in $F_0$ and in~$P$ on the sphere $S^2=\dd
U(O)$ and cut $S^2$ into two disks $S^2_-$ and $S^2_+$ along the former link
(the equator). The mutual position of $\lk_{F_0}O$ and $\lk_PO$ may be
considered as a real analogue for the Milnor fiber of an isolated singularity
on a complex hypersurface, see~\cite{17}. Figure~17 represents the mutual
position of these links on~$S^2$ for minimax and saddle points on a
two-di\-men\-si\-onal component of~$P$ and for the four cases of Fig.~16, where
$O$ is a minimax point of the restriction of~$p$ to an edge of~$P$. Note that
there are, up to isotopy and symmetries, exactly four different embeddings of a
circle with a diameter into a sphere with an equator, provided that any edge of
the graph intersects the equator at most twice and both vertices lie in the
same hemisphere.

The same argument works if a singular point is a vertex~$V$ of a simple
polyhedron~$P$. In this case, $\lk_P(V)$ is a circle with three radii; denote
this graph by~$\Delta$. Any edge of $\Delta$ is the link of~$V$ in one of the
six two-di\-men\-si\-onal components of $(P\sm SP)\cap U(V)$, and we assumed
that $V$ is not a critical point of the restriction of~$p$ to the closures of
these components. This implies that any edge of $\Delta$ intersects the equator
at most twice, because it is close to an arc of a great circle obtained as the
intersection of $S^2=\dd U(V)$ and the tangent plane to a component~$\sigma$
at~$V$. Since the differential $dp|_{\ol e}$ does not vanish at~$V$ for any
edge $e$ of~$P$, the vertices of~$\Delta$ lie in either $S^2_+$ or~$S^2_-$ but
not at the equator $\lk_{F_0}(V)$.

Without loss of generality, we can assume that $P\cap U(V)$ consists of plane
pieces. Now it follows that any edge of $\Delta$ connecting vertices in
different hemispheres intersects the equator exactly once.  An edge with both
endpoints in the same hemisphere can have either none or two intersection
points with the equator; in the latter case, these two points are opposite. An
edge is said to be {\it long\/} if it intersects the equator twice. We have to
consider three cases.

\midinsert

\epsfxsize=362pt

\centerline{\epsffile{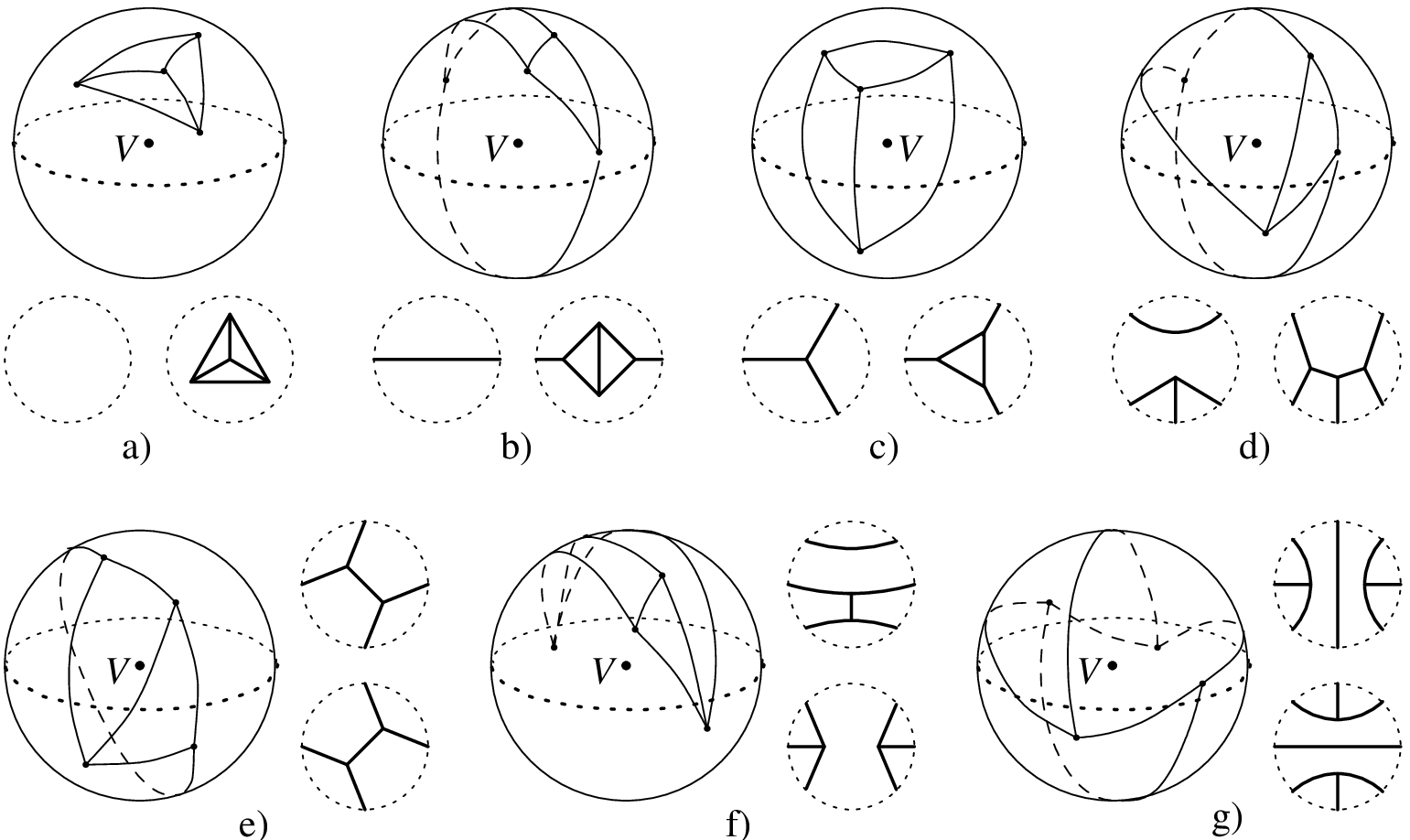}}

\botcaption{Figure 18} Transformation of $K$ near a vertex~$V$  \endcaption

\endinsert

Case~1. All four vertices of~$\Delta$ lie in one hemisphere. Then any long edge
contains a diameter of the other hemisphere. Since any two diameters intersect
one another, there is at most one long edge. Thus in Case~1 there are, up to
isotopy, only two possibilities to draw~$\Delta$ on a sphere with the equator,
see Fig.~18\,a,b.

Figure~18 also shows the difference between $K_-$ and $K_+$.

Case~2. Three vertices of~$\Delta$ lie in one hemisphere and the last one in
the other. Again, there is at most one long edge and thus only two nonisotopic
embeddings of~$\Delta$, see Fig.~18\,c,d.

Case~3. There are two vertices in upper hemisphere and two other in the lower
one. The maximal number of long edges is two: at most one in each hemisphere.
If there are no long edges, we have the situation of Fig.~18\,e; if there is
one long edge, the picture is like Fig.~18\,f; finally, if there are two long
edges, we obtain the situation of Fig.~18\,g.

Note that all results of Section~3.1 apply to the level sets of Morse functions
on arbitrary simple polyhedra, not only to simple spines of~$M(\A)$.

\subhead 3.2. $\th$-curves in the fibers	\endsubhead

Recall that a spine $P$ of a closed 3-manifold~$M$ intersects any loop that is
nontrivial in~$\pi_1(M)$.

\proclaim{Lemma~7} For any $t$, the graph $K_t=P\cap F_t$ intersects any loop
representing a nonzero element of $\pi_1(F_t)=\Z^2$.		\endproclaim

\demo{Proof} Consider the exact sequence
	$$\ldots\to\pi_2(S^1)\to\pi_1(T^2)\overset{i_*}\to\longrightarrow
						\pi_1(M^3)\to\ldots\ $$
of the fibering $p\:M^3\overset{T^2}\to\longrightarrow S^1$. Since $\pi_2(S^1)=
0$, it follows that $i_*$ is a monomorphism, that is, any nontrivial loop in $T
^2$ is nontrivial in $M^3$, too. Hence any spine of~$M$ intersects this loop,
and the Lemma follows.						\qed\enddemo

We assumed that $p$ is an $S^1$-valued Morse function on~$P$. In this case, all
critical points of~$p$ are isolated, and all but finite number of values of~$p$
are regular. Let $F$ be a fiber corresponding to a regular value of~$p$. Then
the intersection $F\cap P$ is a trivalent graph $K\subset F$.

\proclaim{Lemma~8} Suppose that a trivalent graph $K\subset T^2$ intersects any
nontrivial loop in~$T^2$. Then $K$ contains a subgraph that is a $\th$-curve.
								\endproclaim

\demo{Proof} Let $K_1,\dots,K_n$ be connected components of~$K$. If a component
$K_i$ contains no cycles nontrivial in~$T^2$, then there is a disk $D^2_i
\subset T^2$ such that $K\cap D^2_i=K_i$. If there are several components $K_i$
without nontrivial cycles, then there exists a disjoint union of disks~$U$ such
that $K\cap U$ coincides with the set of all connected components of~$K$
containing no nontrivial cycles. Put $K'=K\sm U$. All connected components
of~$K'$ (if any) contain nontrivial cycles.

Any cycle in $\pi_1(T^2)$ is homotopic to a cycle contained in $T^2\sm U$. This
means that any nontrivial cycle intersects~$K'$, and thus $K'\ne\pusto$. Since
any component $K_i'$ of~$K'$ contains nontrivial cycles, we can choose some
cycles $\gamma_i\subset K_i'$ represented by simple closed curves. If there are
several connected components of~$K'$, the cycles $\gamma_i$ do not intersect
one another and thus are homotopic. In this case, among connected components of
$T^2\sm K_1'$ there is an annulus containing other components of~$K'$. It
contains a cycle~$\gamma$ (homotopic to all the~$\gamma_i$) going along a
boundary circle of the annulus. The cycle~$\gamma$ is nontrivial and does not
intersect~$K$. Since this is impossible, the graph~$K'$ is connected.

Let $U_1,\dots,U_k$ be the connected components of $T^2\sm K'$. Every $U_i$ is
an orientable surface with boundary. It contains no closed curves nontrivial
in~$U_i$. Otherwise, such a curve would be either trivial or nontrivial in~$T^
2$. In the former case, it would split the torus into two disjoint parts that
contain different connected components of~$K'$, which is impossible since the
graph~$K'$ is connected. The latter case is impossible, too, because any cycle
nontrivial in~$T^2$ intersects~$K'$. A surface with boundary containing no
nontrivial cycles is a disk. Thus $T^2\sm K'$ is a collection of $s$
2-di\-men\-si\-onal cells, and $K'$ defines a cell decomposition of~$T^2$.

Let $v$, respectively, $e$, be the number of vertices, respectively, edges
of~$K'$. Note that $e=\frac32v$, because $K'$ is a trivalent graph. Also, we
have $v-e+s=\chi(T^2)=0$, which implies that $v=2s$ and $e=3s$. If the number
$s$ is greater than~1, it can be decreased by deleting an edge separating two
different cells. When $s=1$, the graph $K'$ is a $\th$-curve.	\qed\enddemo

\proclaim{Theorem~13} Let $F$ be a nonsingular fiber of the fibering $p\:M^3
\overset{T^2}\to\longrightarrow S^1$, and $K=P\cap F$. Then $K$ contains a
subgraph $L$ that is a $\th$-curve. The number of pairwise nonisotopic \thcs
contained in~$K$ is finite.					\endproclaim

\demo{Proof} Combine Lemmas 7 and 8 and note that the number of all subgraphs
of~$K$ is finite.						\qed\enddemo

Let $F_-$ and $F_+$ be two close nonsingular fibers and $L_-\subset K_-$ be a
\thc in~$F_-$. If the interval $(t_-,t_+)$ contains no singular values of~$p|_
P$, then the graphs $K_-$ and $K_+$ are isotopic, and $K_+$ contains a \thc
$L_+$ isotopic to $L_-$. If there is a singular value $t_0\in(t_-,t_+)$, then
$K_+$ may not contain a \thc isotopic to $L_-$, see, for example, Fig.~19,
where a saddle Morse transformation occurs in a 2-com\-po\-nent of~$P$ as~$t=t_
0$. By Theorem~13, $K_+$ still contains $\th$-curves. It is important that some
of them are not too distant from $L_-$ in the graph~$\Gamma$, that is, can be
obtained from $L_-$ by a small number of flips.

\midinsert

\epsfxsize=360pt

\centerline{\epsffile{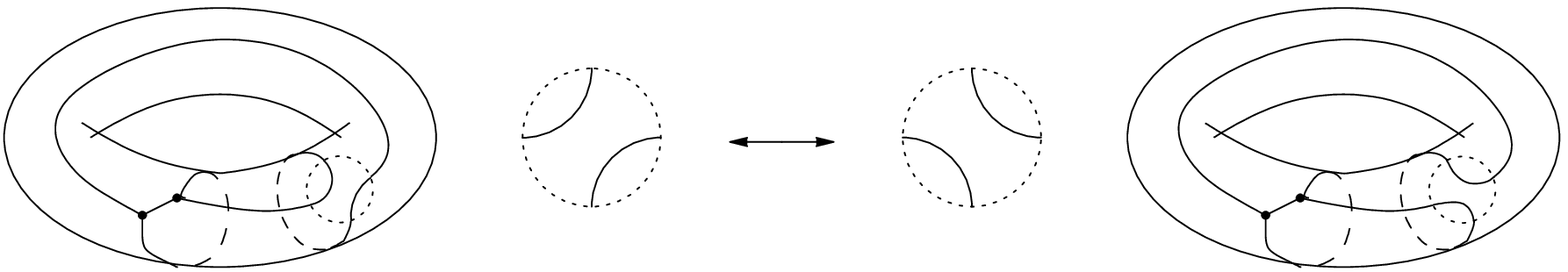}}

\botcaption{Figure 19} Transformation of $L$ induced by a saddle point
\endcaption

\endinsert

\proclaim{Theorem~14} Let both trivalent graphs $K_-,K_+\subset T^2$ contain
\thcs and differ by one of the transformations shown on Figs.~\rom{15},
\rom{16}, and~\rom{18}. Then for any \thc $L_-\subset K_-$ there exists a \thc
$L_+\subset K_+$ such that $d(W_+,W_-)\le5$, where $W_\pm$ is the hexagon
corresponding to~$L_\pm$ and $d$ is the distance in~$\Gamma$, see
Section~\rom{2.1}.						\endproclaim

For a graph $K$, we can consider the set $w(K)\subset\Gamma$ of the hexagons
in~$\Bbb Z^2$ that correspond to all \thcs contained in~$K$. Theorem~14
proclaims that $d(w(K_-),\amb w(K_+))\le5$, where $d(X,Y)$ is the distance from
a subspace~$X$ to a subspace~$Y$ of a metric space~$\Gamma$ defined by the
formula $d(X,Y)=\max\limits_{x\in X}\min\limits_{y\in Y}d(x,y)$; note that $d(X
,Y)$ need not be symmetric.

\demo{Proof} First, suppose that the whole graph $\lk_PV$ lies in one
hemisphere, where $V\in F_0$ is a singular point of $p|_P$ or a vertex. This
happens in the situations of Fig.~18\,a and of the leftmost pictures of Figs\.
15 and~16. Then $U(V)$ contains only an isolated connected component of $K_+$
(or~$K_-$), which has nothing to do with \thcs in~$K_\pm$. So the
transformation of~$K$ arising from a singularity at~$V$ does not affect~$L_-$,
and the graph $K_+$ contains a \thc $L_+$ isotopic to~$L_-$. Of course, in this
case $d(W_+,W_-)=0$.

In the cases represented by Fig.~18\,b and by the second picture of
Fig.~16, it is easy to see that the part of any \thc lying inside of any dotted
circle is a part of an edge of the $\th$-curve (if nonempty), and the part of
the \thc $L_-$ lying outside of the dotted circle can be augmented with a path
in~$K_+$ lying inside of the dotted circle and homotopic to the path in~$K_-$
involved in the $\th$-curve~$L_-$. Similar argument works for the case of
Fig.~18\,c, where the part of~$L_-$ lying inside of the dotted circle may be
a tripod (a neighborhood of a trivalent vertex) or a part of an edge of~$L_-$
or the empty set. So, $d(W_+,W_-)=0$ in all the cases where $\lk_PV$ intersects
the equator in less than four points.

The following Lemma is necessary to deal with the seven remaining cases
(Fig.~18\,d--g, two last pictures of Fig.~16, and the saddle transformation
shown on Fig.~15, right).
\enddemo

\proclaim{Lemma~9} Let $K_-,K_+\subset T^2$ be trivalent graphs and $L_-\subset
K_-$ be a $\th$-curve. Suppose that $K_-$ differs from $K_+$ by one edge only
and $K_+$ contains a $\th$-curve, too.

\rom{1)} If $K_+$ is obtained from $K_-$ by adding one extra edge, then it
contains a \thc $L_+$ such that $d(W_+,W_-)=0$.

\rom{2)} If $K_+$ is obtained from $K_-$ by deleting an edge $e$, then it
contains a \thc $L_+$ such that $d(W_+,W_-)\le1$.		\endproclaim

\demo{Proof} The first statement is obvious, since we can put $L_+=L_-$. In the
second case, we also can put $L_+=L_-$ unless $e\subset L_-$. Suppose that the
edge $e$ is a part of an edge $l$ of~$L_-$. Cut the torus $T^2$ along $L_-$
into a hexagonal 2-cell~$H$. The boundary $\dd H$ contains two arcs $e_1\subset
l_1$, $e_2\subset l_2$ arising from~$e$, and two arcs $f_1=F\!AB$, $f_2=CDE$
complementary to the arcs $l_1=FE$, $l_2=BC$ that arise from~$l$, see Fig.~20.

\midinsert

\epsfysize=80pt

\centerline{\epsffile{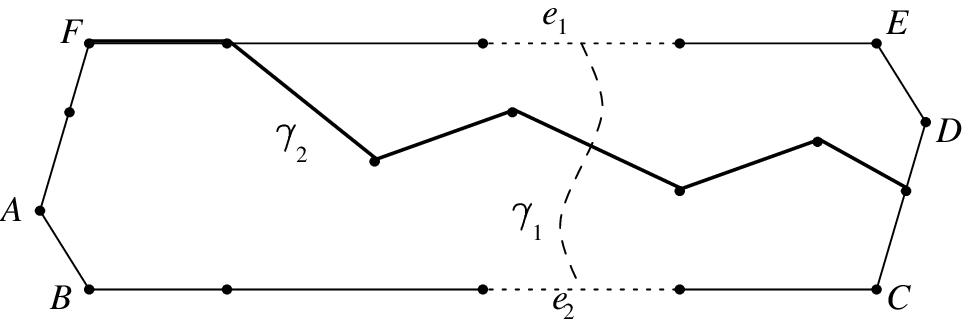}}

\botcaption{Figure 20} The hexagon $H=T^2\sm L_-$		\endcaption

\endinsert

There is the following alternative: either there exists a path $\gamma_1$ in~$H
\sm K_+$ from the midpoint of~$e_1$ to the midpoint of~$e_2$, or there exists a
path~$\gamma_2$ in $K_+$ from a point in~$f_1$ to a point in~$f_2$; this path
$\gamma_2$ may include edges of $K_+$ belonging to~$l\sm e$. The first case
is impossible, because the path $\gamma_1$ yields a nontrivial cycle in $T^2\sm
K_+$; so we have the second case.

Put $L_+=(L_-\sm l)\cup\gamma_2$. Obviously, $L_+$ is a trivalent subgraph
of~$K_+$ with two vertices at the endpoints of~$\gamma_2$. It contains two
edges of~$L_-$ different from~$l$, which form a nontrivial cycle~$\sigma$
homotopic to~$\gamma_1$. Thus $L_+$ intersects any cycle $m\sigma+n\mu\in\pi_1(
T^2)$ with $n\ne0$, where $\mu$ is any cycle such that $\sigma$ and $\mu$
generate~$\pi_1(T^2)$. By construction, $L_+$ intersects any cycle homotopic to
$\sigma$, too. Thus, $L_+$ is a $\th$-curve.

The path~$\gamma_2$ considered above divides~$\dd H$ in two arcs. If the
vertices $A,D$ of $H$, which correspond to two different vertices of~$L_-$ (see
Fig.~20), belong to the same arc, then $L_+$ is isotopic to~$L_-$ and $W_+=W_
-$. If they belong to different arcs, then~$d(W_+,W_-)=1$.  \qed\enddemo

Let us return to the proof of Theorem~14. Consider the case of Fig.~18\,e.
Realize the graph $K_-$ as the union of solid and dotted lines on the middle
picture of Fig.~21. First, add the dashed edge to $K_-$. Then delete the dotted
edge from the resulting graph. This results in~$K_+$. By Lemma~9, the first
step does not affect any \thc in $K_-$, and the graph obtained after the second
step contains a \thc $L_+$ such that $d(W_+,W_-)\le1$.

\midinsert

\epsfysize=60pt

\centerline{\epsffile{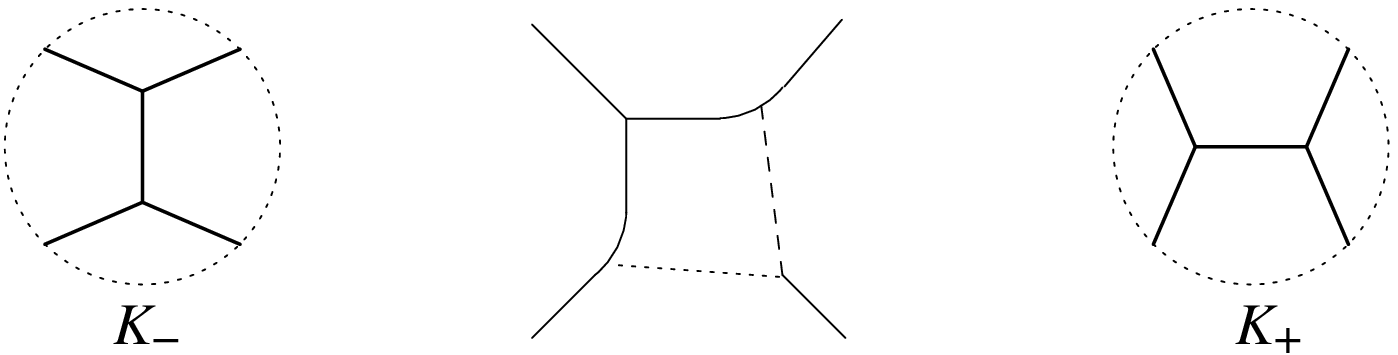}}

\botcaption{Figure 21} Step by step transformation of $K_-$ to $K_+$
\endcaption

\endinsert

The same reasoning applies in all other cases. For instance, for the saddle
transformation (see Fig.~15, right) we need first to add two extra edges that
will be included in~$K_+$, and then delete two edges of $K_+\sm K_-$ from the
graph obtained. Lemma~9 is applicable at each step, because any auxiliary graph
contains either $K_-$ (whenever an edge is added at the previous step) or $K_+$
(whenever an edge is deleted at the next step) and thus contains $\th$-curves.
So for the saddle transformation we can guarantee that $d(W_+,W_-)\le2$ for
some \thc $L_+\subset K_+$. Figure~19 shows that this estimate is exact.

It is easy to see that in all remaining cases the graphs $K_-$ and $K_+$ can be
related by a sequence of $m$ operations of adding an edge followed by $n$
operations of deleting an edge, where $m,n\le5$. In fact, there are sequences
giving the pairs $(m,n)$ equal to $(2,1)$ and $(4,3)$ for the last two pictures
of Fig.~16, $(2,1)$ for the case of Fig.~18\,d, $(3,3)$ and $(5,5)$ for the
cases of Fig.~18\,f and Fig.~18\,g. By Lemma~9, there exists a \thc $L_+\subset
K_+$ such that $d(W_+,W_-)$ does not exceed the number of deletion steps, which
is $n$ for ``left to right'' (or ``bottom to top'') transformations shown on
Figs.~15, 16, and~18, and $m$ for reverse transformations. In any case, this
number is at most~5, which proves the Theorem.				\qed

\remark{Remarks} Theorem~14 gives an upper bound for $d(w(K_+),w(K_-))$. As we
already saw (Fig.~19), this estimate is exact in the case of the saddle
transformation. In fact, there are examples showing that the upper bound
obtained in Theorem~14 is attainable for the third picture of Fig.~16, for the
cases shown on Fig.~18\,d,e, and for the ``top to bottom'' transformation of
Fig.~18\,f. In particular, the distance $d(w(K_+),w(K_-))$ is not symmetric for
some graphs $K_+$ and $K_-$ related by a transformation shown on Fig.~18\,d or
on the third picture of Fig.~16.

The estimate for $d(w(K_+),w(K_-))$ given in Theorem~14 exceeds~3 only in the
case of Fig.~18\,g and for the ``right to left'' transformation from the
rightmost picture of Fig.~16. Even in these cases, no examples are known where
$d(w(K_+),w(K_-))>3$. We believe that such examples do not exist. This would
imply that the upper bound~5 in Theorem~14 can be replaced by~3; the existing
examples show that it cannot be replaced by a number smaller than~3.
\endremark

\subhead 3.3. Proof of Theorem~12				\endsubhead

Let $t\in\R/\Z$ be a parameter on the circle $S^1=\R/\Z$. Without loss of
generality, it can be assumed that the fiber $F_0=p^{-1}(0)$ contains no
singular points of~$p$.

The transversality hypothesis in Theorem~12 means that $p$ has singular points
neither inside of 2-components nor inside of edges of~$P$. Thus the only
singularities of $p|_P$ are the vertices of~$P$; we can assume that the
vertices lie in pairwise different fibers. Let $t_1<t_2<\ldots<t_n$ be their
projections on $S^1$ (that is, the singular values of~$p$); note that $n=c(P)$.
For any $i=1,\dots,n-1$ consider a nonsingular fiber $F_i=F_{(t_i+t_{i+1})/2}$.
Also put $F_n=F_0$.

By virtue of Theorem~13, the graph $K_0=P\cap F_0$ contains a finite number $N>
0$ of $\th$-curves. Let $L_0$ be one of them. By Theorem~14, we can construct
step by step a sequence $L_i\subset K_i=P\cap F_i$, $i=1,\dots,n$, of \thcs
such that $d(W_i,W_{i-1})\le5$. This implies the inequality $d(W_n,W_0)\le5n$.

Recall that $L_n$ is one of the $N$ \thcs contained in~$K_0$. First suppose
that $L_n$ coincides with~$L_0$. For the corresponding hexagons $W_0$ and $W_n$
this means that $W_n=AW_0$; here $A$ is the matrix of the monodromy
operator~$\A$ in the basis where $W_0$ is the standard hexagon. Consequently,
$d(W_0,W_n)=c(A)\ge c(\A)$. Combining this with the inequality $d(W_n,W_0)\le5
n$, we get $c(\A)\le5n$ or, equivalently, $c(P)\ge c(\A)/5$.

Note that there are at least two {\it loose\/} vertices of~$P$, that is,
vertices where the graph~$K$ undergoes the transformation shown on Fig.~18\,c.
One of them corresponds to the minimum point of $\tilde p|_{\ol D}$ and another
to the maximum point, where $\tilde p\:\wt M(\A)\to\R$ is the function that
covers~$p$ under the universal covering $\wt M(\A)\to M(\A)$ and $D$ is one of
the preimages of the 3-cell $M(\A)\sm P$ under this covering. Since loose
vertices do not affect \thcs and corresponding hexagons, we have $d(W_n,W_0)\le
5(n-k)$, where $k\ge2$ is the number of the loose vertices of~$P$. This proves
the inequality $c(P)\ge\frac15c(\A)+k\ge\frac15c(\A)+2$ in the case $L_n=L_0$.

Now suppose that $L_n$ does not coincide with~$L_0$. Consider a sequence $L_{n+
1},\amb L_{n+2},\amb\dots,\amb L_{Nn}$, where $L_i\subset K_{i'}=P\cap F_{i'}$
with $i'\equiv i\mod n$, and $d(W_i,W_{i-1})\le5$ for all $i=1,\dots,Nn$; we
also assume that $W_i=W_{i+1}$ whenever the vertex lying between $F_{i'}$ and
$F_{i'+1}$ is a loose vertex. Since the graph $K_0$ contains only $N$ different
$\th$-curves, we have $L_{sn}=L_{tn}$ for some $s,t$ such that $0\le s<t\le N$.
Then, as above, we have $c(\A^{t-s})\le d(W_{sn},W_{tn})\le5(t-s)(n-k)$,
because $W_{tn}=\A^{t-s}W_{sn}$ and there are $k\ge2$ loose vertices. By
Theorem~8, we have $c(\A^{t-s})=(t-s)c(\A)$ whenever monodromy operator~$\A$ is
non-periodic. Thus $5(t-s)(c(P)-k)\ge(t-s)c(\A)$, i.e., $c(P)\ge\frac15c(\A)+k
\ge\frac15c(\A)+2$. The proof ot Theorem~12 is complete.		\qed

\remark{Remark} It was shown in~\cite{21} that a spine~$P$ can be deformed so
that the only transformations of the level sets are two Morse transformations
(see Fig.~15), the transformations shown on the second picture on Fig.~16, and
flips, see Fig.~18\,d. (In fact, this statement is proved in~\cite{21} for
functions with values in the interval $[0,1]$, but the proof applies to $S^
1$-valued functions, too.) \, Note that a flip in $K$ induces at most one flip
in~$L$ and a saddle transformation induces at most two flips, while two other
transformations do not require flips at all. This yields the estimate $c(P)\ge
c(\A)+2-2s$ (where $s$ is the number of saddle points of~$p$), which holds for
arbitrary spine of~$M(\A)$, implying $c(M(\A))\ge c(\A)+2-2s$. However, it is
unlikely that there exists any lower bound for~$s$.		\endremark

\head\S4. Main conjectures \endhead

\subhead 4.1. Triangulations of $T^2$ coming from spines \endsubhead

Let us recall the idea of the proof of Theorem~12. Given a spine $P$ of $M^3=M^
3(\A)$, suppose that we can assign some object $\sigma_t$ to each fiber $F_t$
containing no vertices of~$P$ in such a way that
\roster
\item"1)" $\sigma_t$ is an element of some metric space~$(\Sigma,d)$;
\item"2)" there is an action of the group $\SL(2,\Z)$ on the set~$\Sigma$, and
	monodromy $\A$ takes $\sigma_t$ to $\sigma_{t+2\pi}=\A\sigma_t$;
\item"3)" $d(\sigma_t,\A\sigma_t)\ge c(\A)$;
\item"4)" as $t$ varies, $\sigma_t$ remains unchanged until $F_t$ encounters a
	vertex $V\in P$;
\item"5)" if there is exactly one vertex of $P$ between $F_{t_+}$ and $F_{t_-
	}$ (we can assume that vertices of~$P$ lie in pairwise different
	fibers), then $d(\sigma_{t_+},\sigma_{t_-})\le1$.
\endroster
Then the number of vertices of~$P$ is at least~$c(\A)$. This can be followed by
an argument that shows that this number is in fact at least $c(\A)+k$, where
$k$ is some positive integer smaller than~6.

In the proof of Theorem~12, $\sigma_t$ is an isotopy class of \thcs in~$F_t$,
$\Sigma$ is the vertex set of the trivalent graph~$\Gamma$ described in
Section~2.1, $d$ is the distance in~$\Gamma$, and~$k=2$. To preserve
property~4, we had to restrict ourselves to spines transversal to the fibers.
Still, we got $d(\sigma_{t_+},\sigma_{t_-})\le5$ instead of property~5; this
leads to an annoying factor of~$1/5$ in Theorem~12.

In this section we present another construction for~$\sigma_t$; namely, it will
be some class of 2-chains representing the fundamental class~$[F_t]$. This
approach applies to all special spines~$P$ of~$M(\A)$ (recall that, by
Theorem~2, any minimal spine of $M(\A)$ is a special one), not only to those
transversal to the fibers, and properties 1, 2, 4, and~5 hold, while property~3
remains unproved.

Given a special spine~$P$ of~$M^3$, by $P'$ denote a triangulation of~$M^3$
dual to~$P$. Its tetrahedra correspond to the vertices of~$P$ (see Fig.~1\,f),
triangles correspond to the edges of the singular graph $SP$ (that is, to the
triple lines) of the spine~$P$, edges of~$P'$ are dual to the 2-cells of~$P$,
and the only vertex of~$P'$ is located somewhere inside of the 3-cell~$M\sm P$.
Here and below we abuse the word ``triangulation'': the intersection of two
simplices (if nonempty) is a subset of the sets of their faces of smaller
dimensions, but the cardinality of this subset may exceed~1.

Let $f\:\R^1\to S^1$ be the standard covering defined by the rule $f(t)=t\mod2
\pi$. By $\tilde f\:\tM^3\to M^3$ denote the covering such that $\tilde f_*(\pi
_1(\tM^3))=i_*(\pi_1(F))$, where $i\:F\to M^3$ is the embedding map of some
fiber~$F$ into~$M^3$. Then a projection $\tilde p\:\tM\to\R^1$ can be defined
by the relation $f\!\circ\tilde p=p\circ\!\tilde f$. By $F_t\subset\tM$ denote
the fiber $\tilde p^{-1}(t)$, where $t\in\R^1$. Put $\wt P=\tilde f^{-1}(P)$
and $\wt P'=\tilde f^{-1}(P')$. Note that $\wt P$ is a spine of $\tM$ punctured
infinitely many times (at all preimages of the vertex~$V'$ of~$P'$) and $\wt P
'$ is a triangulation of~$\tM$. Furthermore, $\tM$ is nothing but $T^2\x\R^1$
equipped with the deck transformation $(x,t)\mapsto(\A x,t+2\pi)$, where $x\in
T^2$ and $t\in\R^1$, and both $\wt P$ and $\wt P'$ are invariant under the deck
transformation. Fix a cartesian product structure on $\tM$ and define the
forgetting projection $j\:\tM\to\T^2$ by $j(x,t)=x$.

By $M_0$, respectively, $\tM_0$, denote the complement to the vertices of~$P$
in~$M$, respectively, to the vertices of~$\wt P$ in~$\tM$. It is easy to see
that there is a strict deformation retraction $r\:M_0\to\sk_2P'$, which is
covered by a strict deformation retraction $\tilde r\:\tM_0\to\sk_2\wt P'$. If
a fiber $F_t$ contains no vertices of~$\wt P$, the image $\tilde r(F_t)\subset
\sk_2\wt P'$ defines a 2-cycle $c^2_t\in Z_2(\wt P')$. The group $Z_2(\wt P')$
is generated by the 2-cells of~$P'$, so we have $c^2_t=\sum\a_k(t)\tau^2_k$,
where $\{\tau^2_k\}$ is the set of the 2-cells of~$P'$.

The coefficients~$\a_k(t)$ are equal to the intersection indices of~$F_t$ with
the oriented edges $e_k$ of~$\wt P$ corresponding to the (oriented) 2-cells
$\tau^2_k$. Let $A_k$ and~$B_k$ be the endpoints of an edge~$e_k$ oriented
from $A_k$ to~$B_k$. If $A_k$ and $B_k$ lie both below or both above~$F_t$
(that is, if $(\tilde p(A_k)-t)(\tilde p(B_k)-t)>0$), then $\a_k(t)=0$;
otherwise, $\a_k(t)=\sign(\tilde p(B_k)-\tilde p(A_k))=\pm1$. Obviously, only
finite number of the coefficients can be nonzero, and property~4 holds for
all~$\a_k(t)$, whence for~$c^2_t$ as well.

Let us consider a 2-chain $j_*(c^2_t)$. By construction, it is a cycle and its
homology class is $[T^2]$ (because $j_*(c^2_t)=j_*(F_t)$). All triangles that
constitute $j_*(c^2_t)$ have all their vertices mapped to $j(V')$ (where $V'$
is the vertex of~$P'$). Let $\sigma_t$ be the linearization of $j_*(c^2_t)$,
that is, the 2-chain obtained from $j_*(c^2_t)$ by replacing the characteristic
mappings of all triangles by homotopic (with fixed vertices) linear mappings; a
mapping of a triangle to the torus is linear if it is a linear mapping to the
plane followed by the projection~$\R^2\to T^2$. Roughly speaking, we define
$\Sigma$ as the set of all linear 2-cycles in~$T^2$ that are homological to~$[T
^2]$ and have the only vertex, which is placed at~$j(V')$. Accurately speaking,
any element of~$\Sigma$ is a (homological to~$[T^2]$) cycle $\sum\beta_k\Delta_
k$, where the coefficients $\beta_k$ are integers, all but finite number of
them are zero, and $\{\Delta_k\}$ is the set of projections (from~$\R^2$ to~$T^
2$) of all different oriented triangles with vertices in~$\Z^2$; the triangles
are equal if they come from the same edge $e_k$ of~$\wt P'$, so even
geometrically coinciding triangles may be different in this sense. This
difference will be expoited later, but unless otherwise stated, we identify all
coinciding triangles. Degenerated triangles (those with three vertices along a
line) are allowed, and their orientation is, as usually, a cyclic ordering of
their vertices. The metric on~$\Sigma$ will be constructed later.

Obviously, we have $\sigma_t\in\Sigma$. Further, the group $\SL(2,\Z)$ acts
on~$\Sigma$ in a natural way, and $\sigma_{t+2\pi}=\A\sigma_t$, where $\A$ is
the monodromy of the fibration $p\:M^3\to S^1$. Thus properties~2 and~4, see
above, are satisfied; further, property~1 is satisfied for any choice of
metric~$d$.

To find the difference between $\sigma_{t_+}$ and $\sigma_{t_-}$ (where the
interval $(t_-,t_+)$ contains the projection $\tilde p(V)$ of exactly one
vertex $V$ of~$\wt P$), consider the difference between $F_{t_+}$ and~$F_{t_-
}$. It is easy to see that $F_{t_+}$ can be obtained as a connected sum of $F_{
t_-}$ and a small two-dimensional sphere $S^2(V)$ centered at~$V$. Thus $\sigma
_{t_+}$ is obtained from~$\sigma_{t_-}$ by adding $j_*(\tilde r(S^2(V)))$, that
is, by one of the two-dimensional Pachner moves, see~\cite{20} and \S2.3~above.
Indeed, there are four edges of~$\wt P$ emanating from~$V$. Suppose that they
are all different, that is, they are not loops. Then~$\tilde r(S^2(V))$ is the
boundary of a tetrahedron triangulated into four triangles~$\tau^2_k$, from
which 0 to 4 can annihilate with triangles that constitute~$\tilde r(F_{t_-})$.
So~$\sigma_{t_+}$ is obtained from~$\sigma_{t_-}$ in one of the following five
ways:
\roster
\item"0)" by adding the projection of the boundary of some tetrahedron, that
	  is, by adding triangles $BCD$, $CAD$, $ABD$, and~$BAC$ (mind the
	  orientations!), or
\item"1)" by replacing a triangle $ABC$ of~$\sigma_{t_-}$ by three triangles
	  $ABD$, $BCD$, and $CAD$, or
\item"2)" by replacing two triangles $ABC$ and $ACD$ of~$\sigma_{t_-}$ (the
	  segment $AC$ contributes to both their boundaries, but with opposite
	  signs) by the triangles $ABD$ and~$BDC$, or
\item"3)" by replacing three triangles $ABD$, $BCD$, and~$CAD$ by a single
	  triangle~$ABC$, or
\item"4)" by erasing four triangles $BCD$, $CAD$, $ABD$, and~$BAC$.
\endroster
Of course, moves 3 and 4 may not be applicable to arbitrary~$\sigma_t\in
\Sigma$.

If one of the edges of $\tilde P$ incident to~$V$ is a loop, then two of the
triangles that form $\tilde r(S^2(V))$ annihilate with one another, and it is
easy to see that $\sigma_{t_+}=\sigma_{t_-}$.

Consider a graph with vertices at the elements of~$\Sigma$ and edges
corresponding to the moves described above. Let $d$ be the distance function on
this graph. This completes the construction of the metric space $(\Sigma,d)$.
Property~5 now holds by construction. Recall that properties 1, 2, and~4 hold
as well. Property~3 is a conjecture. If it holds, then the inequality $c(M(\A))
\ge c(\A)$ follows immediately.

This conjecture is supported by the following observation. For any triangle
of~$\sigma_t$, draw three segments connecting its baricenter with the midpoints
of its sides. All the tripods obtained in this way form a graph~$K$. In some
simple cases (for example, for spines constructed in Section~2.3), $K$~is a
trivalent graph embedded into~$T^2$. The Pachner moves 0 and~4 affect $K$
according to Fig.~18\,a (``left to right'' transformation of Fig.~18\,a for
move~0 and ``right to left'' for move~4); moves 1 and 3 affect $K$ according to
Fig.~18\,c; finally, moves~2 correspond to flips, see Fig.~18\,e. It is
possible now to choose a \thc in~$K$ and apply the argument of~\S3, taking into
account that moves of Figs.~18\,a and~18\,c do not affect the isotopy classes
of \thcs and any flip in the graph~$K$ implies at most one flip of any \thc in
it. However, in the general case the graph~$K$ neither is embedded into the
torus nor is trivalent, because there may be 4, 6 etc\. triangles having the
same edge in common, even for spines transversal to the fibers.

There is a more geometrical reformulation of this conjecture. For a chain $c=
\sum_ir_i\tau_i$, where the sum is finite, $r_i\in\Z$, and the $\tau_i$ are
singular simplices, put $\|c\|=\sum_i|r_i|$, and consider the minimum value~$l^
1(\A)$ of~$\|c\|$ over all 3-chains in $T^2\x I$, $I=[0,2\pi]$, representing
the fundamental class $[T^2\x I]$ and satisfying the condition $\sigma_{2\pi}=
\A\sigma_0$, where 2-chains $\sigma_{2\pi}$ and~$\sigma_0$ are the
intersections of the chain~$c$ with $T^2\x\{2\pi\}$ and $T^2\x\{0\}$, and $\A(x
,0)=(\A x,2\pi)$ for $(x,0)\in T^2\x\{0\}$. It is easy to see that $l^1(\A)\le
c(\A)$. The conjecture can be formulated as follows: $l^1(\A)=c(\A)$. The
definition of $l^1(\A)$ is similar to Gromov's definition of the simplicial
$l^1$-norm (with some boundary conditions imposed), see~\cite{10}. However,
with the original definition by Gromov, we have $\|M^3(\A)\|_{l^1}=0$. Also see
the concluding remarks in~\cite{26}.

In fact, the estimate $c(M(\A))\ge c(\A)$ can be deduced from Conjecture~3 (see
below), which is weaker than the conjecture above. Recall that the triangles
of~$\sigma_t$ are ``marked'' by the corresponding edges of~$\wt P$. It may
happen that two different edges contribute two equal but oppositely oriented
triangles to~$\sigma_t$. Let us prohibit cancellation of the triangles in these
cases. Then we get the set~$\Sigma'$ of 2-cycles homological to~$[T^2]$ such
that any triangle is marked by two nonnegative integers, say, $(m,n)$, which
means that $\sigma'$ contains $m$ positively oriented copies of this triangle
and $n$ negatively oriented copies of it. A natural mapping $s\:\Sigma'\to
\Sigma$ is defined by replacing a mark $(m,n)$ by the coefficient $m-n$. To
define a distance $d'$ on~$\Sigma'$, we follow the construction of the distance
function~$d$; however, we prohibit moves 0 and~4; only moves 1, 2, and~3 are
allowed. Obviously, the inequality $d'(\sigma_1,\sigma_2)\ge d(s(\sigma_1),s(
\sigma_2))$ holds. Note also that $s(\A\sigma)=\A s(\sigma)$, where $\sigma\in
\Sigma'$ and $\A\in\SL(2,\Z)$. So the following statement is weaker than the
conjecture above (which states that property~3 holds).

\proclaim{Conjecture~3} For all $\sigma\in\Sigma'$ and $\A\in\SL(2,\Z)$, the
inequality $d'(\A\sigma,\sigma)\ge c(\A)$ holds.		\endproclaim

However, this weaker conjecture still implies the inequality $c(M(\A))\ge c(\A
)$.

\proclaim{Theorem~15} Conjecture~\rom3 implies the estimate $c(M(\A))\ge
c(\A)$.								\endproclaim

\demo{Proof} Let us say that a vertex $V$ of $\wt P$ is {\it maximal},
respectively, {\it minimal}, if the endpoints of all four edges going from~$V$
are not above, respectively, not below~$V$ (for $A,B\in\tM$, we say that $A$
lies above $B$ if $\tilde p(A)>\tilde p(B)$; further, we assume that $\tilde p(
A)\ne\tilde p(B)$ whenever vertices $A$ and~$B$ are different). Maximal and
minimal vertices are also called {\it critical}. If a vertex $V\in\wt P$ is
maximal (respectively, minimal, critical), then its projection $\tilde f(V)\in
P$ is said to be a {\it maximal\/} (respectively, {\it minimal}, {\it
critical\/}) vertex of~$P$. Obviously, maximal (minimal, critical) vertices
of~$P$ are well-defined, because for any vertex $V\in P$ all its preimages
$\tilde f^{-1}(V)$ are or are not critical (maximal, minimal) vertices of~$\wt
P$ simultaneously.

First, suppose that there are no critical vertices. Then any vertex of~$\wt P$
induces a change of $\sigma_t\in\Sigma'$ by one of the moves 1, 2, 3.
(However, it can happen that the chain $s(\sigma_t)\in\Sigma$ undergoes the
move 0 or~4. That is why we had to introduce $\Sigma'$ instead of~$\Sigma$.)
\,It follows that in this case $d(s(\sigma_t),s(\sigma_{t+2\pi}))=d'(\sigma_t,
\sigma_{t+2\pi})$, so Conjecture~3 implies property~3 (see the beginning
of~\S4), and the spine~$P$ contains at least $c(\A)$ vertices.

Suppose that $\wt P$ contains critical vertices, but there is an isotopy of the
embedding of $P$ in~$M^3$ such that the deformed spine $P_1$ yields the
deformed covering spine $\wt P_1$ with neither minimal no maximal vertices.
Then the number of the vertices of~$P_1$ is at least~$c(\A)$ (provided that
Conjecture~3 is true); obviously, $P$ and $P_1$ have the same number of
vertices.

The notion of peripheric edge is necessary for the case of arbitrary spines.
\enddemo

\definition{Definition~13} An edge $e$ of the singular graph $SP$ of a
spine~$P$ is said to be {\it peripheric\/} if there exists a vertex $V$ of~$P$
with the following property: if a cycle $\gamma\subset SP$ contains $e$ and $p_
*(\gamma)\ne0$ (here $p$ is the projection $M^3\to S^1$), then every connected
component of $\gamma\sm V$ containing $e$ is a loop $\gamma'$ (with endpoints
at~$V$) such that $p_*(\gamma')=0$; in the other words, the map $p_*$
restricted to the closure of the connected component of $SP\sm V$
containing~$e$ is trivial. Edges that are not peripheric are called {\it
regular}.							\enddefinition

For example, a loop $e$ is peripheric if and only if $p_*(e)=0$. Unlike the
property of a vertex to be critical, the property of an edge to be peripheric
is an isotopy invariant.

\proclaim{Lemma~10} If $P$ contains no peripheric edges, then it is isotopic to
a spine with no critical vertices.				\endproclaim

\demo{Proof} Take any edge $e_1\in P$. Since $e_1$ is a regular edge, there
exists a cycle $\gamma_1\subset SP$ such that $e_1\subset\gamma_1$, $p_*(\gamma
_1)\ne0$, and $\gamma_1$ does not pass twice through any vertex. Then there
exists a spine $P_1$ isotopic to $P$ such that the projection $\tilde p\:\tM\to
\R^1$ restricted to a connected component $\tilde\gamma'_1$ of $\tilde f^{-1}(
\gamma'_1)$, where $\gamma'_1\subset P_1$, takes the sequence of the vertices
of~$\tilde\gamma'_1$ to a strictly monotone sequence $\{a_i\mid i\in\Z\}$ of
real numbers; here $\gamma'_1$ is obtained from~$\gamma_1$ under an isotopy
that takes $P$ to~$P_1$. To construct~$P_1$, fix a monotone sequence $\{a_i\}$
satisfying the condition $a_{i+m}=a_i+2\pi k$, where $m$ is the length
of~$\gamma_1$ and $k\in\Z=\pi_1(S^1)$ is equal to $p_*(\gamma_1)$. Let $\{b_i
\mid i\in\Z\}$ be the sequence of the projections $\tilde p(V_i)$, where the $V
_i$ are consecutive vertices of $\tilde\gamma_1$. The differences $c_k=a_k-b_
k$ form an $m$-periodic sequence. Let $\{\delta_i\mid i=1,\dots,m\}$, be a
family of arcs in~$\tM$ such that the endpoints of $\delta_i$ are $V_i$ and $V'
_i$, where $\tilde p(V'_i)=a_i$ and the $\delta_i$ do not intersect~$\wt{SP}$;
suppose also that the projections $\tilde f(\delta_i)$ do not intersect one
another. It is easy to see that such a family exists, and that it is possible
to get a spine~$P_1$ with the properties described above by an isotopy of the
mapping \,$\id\:M^3\to M^3$, fixed outside of a small collar neighborhood
of~$\bigcup\limits_{i=1}^ms(\delta_i)$.

Now any vertex $V\in\gamma_1$ is not critical. If there remain critical
vertices, take an edge $e_2\in P_1$ incident to one of them and a cycle $\gamma
_2\subset SP_1$ containing~$e_2$ such that $p_*(\gamma_2)\ne0$ and $\gamma_2$
does not pass twice through any vertex. If $\gamma_2\cap\gamma_1=\pusto$,
repeat the construction above. Otherwise, take the connected component $\gamma_
2^\circ$ of $\gamma_2\sm\gamma_1$ that contains~$e_2$, and move the inner
vertices of this component as above, but do not move its endpoints. In the
spine $P_2$ obtained after this step, any vertex $V\in\gamma_1\cup\gamma_2^
\circ$ is not critical. Indeed, any vertex of~$\gamma_2^\circ$ is a noncritical
vertex of~$P_2$ by construction, and any vertex $V\in\gamma_1$ has neighboring
(in~$\gamma_1$) vertices above and below $V$ in~$P_1$; $P_2$ inherits this
property from~$P_1$, because the isotopy converting $P_1$ into~$P_2$ does not
affect~$\gamma_1$.

If there still remain critical vertices, repeat the construction described
above. Then we get a spine $P_3$ isotopic to~$P_2$. The set of the vertices of
$\gamma_1\cup\gamma_2^\circ\cup\gamma_3^\circ$, which certainly are not
critical in $P_3$, is larger than the similar set $\gamma_1\cup\gamma_2^\circ$
in~$P_2$. Hence, in a finite number of steps we obtain a spine~$P_k$ isotopic
to~$P$ and having no critical vertices. 			\qed\enddemo

According to Lemma~10, it follows from Conjecture~3 that any spine of~$M(\A)$
without peripheric edges contains at least~$c(\A)$ vertices.

The degree of any vertex $V\in SP$ equals four. Suppose that there are $k$
peripheric and $4-k$ regular half-edges incident to~$V$ (we consider
half-edges, because the graph~$SP$ may contain loops). Note that the number~$k$
is even, thence is equal to 0, 2 or~4. Indeed, if there is a regular edge $e$
incident to~$V$ (that is, if $k<4$), then there is a cycle $\gamma$ consisting
of regular edges and passing through~$V$, so in this case the number of regular
half-edges incident to~$V$ is at least~2, and we have~$k\le2$. On the other
hand, if $k=1$, consider the connected component of $SP\sm V$ containing the
peripheric edge $e$ incident to~$V$. Let it contain $m$ edges and $n$ vertices
different from~$V$. Then the number of half-edges in this component is equal to
$2m$ and to $4n+1$ simultaneously, which is impossible.

Let $P$ be an arbitrary spine of $M(\A)$. Let us say that a vertex is {\it
regular\/} if it is incident to regular edges only, {\it semiregular\/}
if it is incident to two peripheric and two regular half-edges, and {\it
peripheric\/} if it is incident to peripheric edges only. Following the proof
of Lemma~10, replace $P$ by an isotopic spine $P_1$ such that all regular and
semiregular vertices of~$P_1$ are not critical. For any semiregular vertex~$V$
of~$P_1$, by $b(V)$ denote the connected component of $SP\sm V$ that consists
of peripheric edges only (i.e., lies ``behind~$V$''). Following the proof of
Lemma~10, replace~$P_1$ by a spine~$P_2$ such that $P_2$ is isotopic to~$P_1$,
the isotopy between $P_1$ and~$P_2$ does not affect regular edges, $p(V_i)=p(V
)$ for all semiregular vertices $V\in P_2$ and all peripheric vertices $V_i\in
b(V)$, and $p(e)$ is a ze\-ro-ho\-mo\-to\-pic loop in~$S^1$ for any peripheric
edge~$e$.

Consider 2-cycles $\sigma_t\in\Sigma'$ constructed from the spine~$P_2$. If
there is exactly one regular vertex between $F_{t_+}$ and~$F_{t_-}$, then
$\sigma_{t_+}$ is obtained from~$\sigma_{t_-}$ by one of the moves 1, 2,~3, so
$d'(\sigma_{t_+},\sigma_{t_-})=1$.

Suppose that there is exactly one semiregular vertex~$V\in\wt P_2$ (together
with all peripheric vertices of~$b(V)$) between $F_{t_+}$ and~$F_{t_-}$.
Compare $\sigma_{t_-}$ and~$\sigma_{t_+}$. The former cycle consists of the
triangle corresponding to the regular edge $e_-$ of~$\wt P_2$ incident to~$V$
such that the other endpoint of~$e_-$ lies below~$V$, and of some other
triangles. The latter cycle includes the triangle that corresponds to the
other regular edge~$e_+$ of~$\wt P_2$ incident to~$V$; other triangles that
constitute~$\sigma_{t_+}$ are the same as in~$\sigma_{t_-}$. Note that all
triangles but one of a 2-cycle $\sigma\in\Sigma'$ determine the remaining
triangle~$\Delta$ uniquely: directions of its sides are those of the boundary
of the 2-chain formed by the remaining triangles and thus denote~$\Delta$ up
to dilatation with the coefficient~$\pm1$, and this coefficient is uniquely
determined by the condition that $\sigma$ is homological to~$[T^2]$;
informally, $\Delta$ is the difference between $[T^2]$ and the 2-chain formed
by the other triangles of~$\sigma$. So in this case we have $\sigma_{t_+}=
\sigma_{t_-}$.

We have shown that regular vertices of~$P_2$ shift~$\sigma_t$ by distance~1,
while semiregular and peripheric vertices do not affect $\sigma_t$ at all.
Recall that regularity, semiregularity, and periphericity are preserved under
isotopy. Consequently, Conjecture~3 implies that the number of regular vertices
of an arbitrary spine $P$ of $M(\A)$ is at least~$c(\A)$, so the number of all
vertices of~$P$ is greater than or equal to~$c(\A)$, too. This completes the
proof of Theorem~15.							\qed

Finally, let us provide a sketchy explanation (or, maybe, rather motivation)
for the ``$+5$'' summand in Conjecture~2. Suppose that there is an edge $e\in P
'$ of the triangulation~$P'$ dual to a given spine~$P$ such that $p_*(e)=\pm1$
in the group $\pi_1(S^1)=\Z$. Note that the cartesian product structure on
$\tM$ and the corresponding projection $j\:\tM\to T^2$ are not uniquely
defined; in our case we can choose any element of $\Z^2=\pi_1(T^2)$ to be the
image~$j_*(e)$. Choose $j$ so that $j_*(e)=0$. Cycles $\sigma_t$ can include
degenerated triangles with an edge $j_*(e)$, which are projections of the
triangles of~$\wt P'$ incident to~$e$. However, if vertices $A$ and $B$ of a
triangle~$ABC$ coincide, the triangle~$ABC$ cannot be distinguished from~$BAC$,
which is equal, on the other hand, to~$ABC$ with the opposite orientation, so
such triangles can be ignored in~$\sigma_t$. It can be easily shown that
Pachner moves 1, 2, and~3 involving triangles of this type do not
affect~$\sigma_t$ and thus require $k$ ``additional'' vertices of~$P$, where
$k$ is the number of triangles of~$\wt P'$ incident to~$e$, that is, the number
of sides of the 2-cell of~$\wt P$ dual to~$e\subset\wt P'$. We can assume $P$
to be a minimal spine of~$M(\A)$. In this case $P$ is a special spine, and,
according to~\S1.2, it contains $n+1$ 2-cells, where $n=c(M^3)$ is the number
of its vertices. The total number of the vertices of 2-cells equals~$6n$, so
the average number of sides for a 2-cell of~$P$ is $\frac{6n}{n+1}>5$ (because
$n\ge6$). Thus it is reasonable to expect at least 5 ``spare'' vertices, in
addition to $c(\A)$ vertices that do change~$\sigma_t$.

Of course, it can happen that there are no edges~$e$ such that $p_*(e)=\pm1$.
However, there always exists an edge~$e$ such that $p_*(e)=m>0$. Then the
argument above can be applied to an $m$-fold covering $\tM_m=M(\A^m)$ of~$M(\A
)$. A special spine $P$ of $M(\A)$ (punctured once) is covered by a special
spine $\wt P_m$ of $\tM_m$ with $m$ punctures, and we can destroy $m-1$ 2-cells
of $\wt P_m$ in order to get a special spine of $\tM_m$ punctured only once.
Apart from this, we destroy an $m$\snug th 2-cell of~$\wt P_m$ by the process
described in the previous paragraph. Consequently, we obtain vertices of~$m$
2-cells of~$\wt P_m$ as additional, spare vertices, and $c(\A^m)=mc(\A)$ (see
Theorem~8) vertices of~$\wt P_m$ affecting~$\sigma_t$; obviously, the number of
vertices of~$\wt P_m$ is equal to $m$ times the number of the vertices of~$P$.
Unfortunately, the ``averaging argument'' provides no proof for the claim that,
by destroying~$m$ 2-cells during the above-described process, we can destroy at
least $5m$ vertices (in fact, $4m+1$ would already suffice), even though there
is an essential ambiguity in the choice of 2-cells to be discarded. Because of
pseudominimality of~$P$, all boundary curves of 2-cells are not short,
see~\cite{14} and the beginning of~\S2.3 above, but this does not readily yield
$4m$ of $4m+1$ spare vertices required, since {\sl \`a priori\/} the boundary
curve of a 2-cell can visit some vertices more than once. It is only possible
to prove with this approach that $2m$ vertices are destroyed when $m$ 2-cells
are discarded, and this results in 2 ``additional'' vertices of~$P$, compare
with Theorem~12.

\subhead 4.2\endsubhead {\bf Can hyperbolic geometry help?\/}

Hyperbolic plane~$H^2$ can be a natural receptacle for various metric spaces
$\Sigma$ introduced in~\S 4.1, because $H^2$ is the Teichm\"uller space of~$T^
2$. Being a subgroup of $\SL(2,\R)$, the group $\SL(2,\Z)$ acts on~$H^2$
(modelled on the upper half-plane) by the rule 
	$$\A(z)=\dfrac{az+b}{cz+d},\quad\text{where}\quad\A=\pmatrix a&b\\
		c&d\endpmatrix\in\SL(2,\Z)\quad\text{and}\quad\Im z>0;\tag4$$
the kernel of this action is $\{\pm I\}$, so formula~(4) defines an action of
the modular group $\SL(2,\Z)/\{\pm I\}$ on~$H^2$.

Consider the ideal triangle in~$H^2$ with vertices at $0$, $1$, and~$\infty$.
Take its mirror images in its sides. This gives the triangles $(-1,0,\infty)$,
$(0,1/2,1)$, and~$(1,2,\infty)$, where $(a,b,c)$ denotes the ideal triangle
with vertices $a$, $b$, and~$c$. On the next step, construct the images of the
triangles obtained in the previous step under reflections in their sides that
are not sides of triangles obtained earlier. Continuing this way, we get a
tesselation of~$H^2$ into equal ideal triangles. It is called the {\it Farey
tesselation}. Note that the modular group action~(4) preserves it.

Consider a black-and-white coloring of triangles of the Farey tesselation such
that neighboring triangles are of different colors. The only epimorphism of the
modular group to~$\Z_2$ discussed in Theorem~6 can be defined as follows: $\f(g
)=0$ if $g$ preserves the coloring and $\f(g)=1$ if $g$ reverses it.

The following lemma (which certainly is well known though I did not find a
reference) reveals the connection between the Farey tesselation and the Farey
series.

\proclaim{Lemma~11} Let $r/s$ be the third vertex of a triangle $\Delta_1$
obtained by the reflection in the side $(m/n,p/q)$ of a triangle~$\Delta$
obtained on the previous step. Then $\dfrac rs=\dfrac{m+p}{n+q}$.  \endproclaim

\demo{Proof} The proof is by induction over the number of reflections,
counting from the ``root'' triangle $\Delta_r=(0,1,\infty)$. By the induction
hypothesis, the third vertex of~$\Delta$ is~$\frac{m-p}{n-q}$. A direct
calculation of the cross-ratio of the quadruple $\(\frac mn,\frac{m+p}{n+q},
\frac pq,\frac{m-p}{n-q}\)$ shows that the points $\frac{m-p}{n-q}$ and $\frac{
m+p}{n+q}$ are symmetric in the line $(\frac mn,\frac pq)$.	\qed\enddemo

A line $(m/n,p/q)$, where $m,n,p,q\in\Z$, occurs in the Farey tesselation if
and only if~$mq-np=\pm1$ (we express $\infty$ as~$1/0$). The ``only if''
statement is proved by induction on the number of reflections, as above. To
prove the ``if'' statement, read Chapter~III of~\cite{11} or note that $m/n$
and $p/q$ either are both nonnegative or both nonpositive, since otherwise $|mq
-np|>1$. Hence, the line $(m/n,p/q)$ does not intersect the line~$(0,\infty)$.
Any other line $l$ of the Farey tesselation can be taken to~$(0,\infty)$ by a
sequence of reflections in its neighboring lines followed by a rotation of the
triangle $(0,1,\infty)$. Let this isometry of~$H^2$ take the line $(m/n,p/q)$
to $(m'/n',p'/q')$, then $|m'q'-n'p'|=|mq-np|\ne1$, so the line $(m/n,p/q)$
does not intersect $l$ as well. Then it cannot pass through an interior point
of any triangle of the tesselation, and the only chance for it to be contained
somewhere in~$H^2$ is to be one of the lines of the tesselation. For a
constructive proof of the ``if'' statement, the Euclid algorithm should be
applied to one of the pairs $(m+p,n+q)$ and $(m-p,n-q)$, compare with the proof
of Theorem~2 in~\cite2.

\midinsert

\epsfxsize=330pt

\centerline{\epsffile{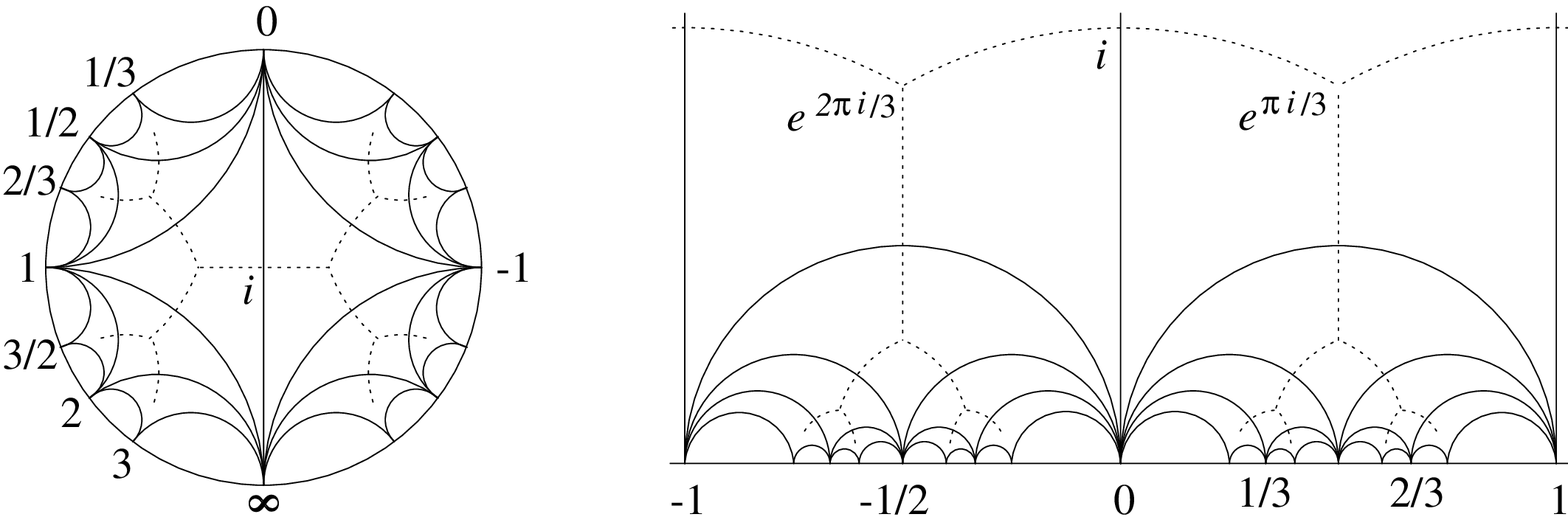}}

\botcaption{Figure 22} Farey tesselation of hyperbolic plane	\endcaption

\endinsert

Figure~22 represents the Farey tesselation of~$H^2$ for the circle and
half-plane models of hyperbolic plane. For more information on the Farey
tesselation see~\cite{23}.

Recall that isotopy classes of \thcs correspond to admissible hexagons,
see~\cite2 or~\S2.1 above. There is a natural correspondence between admissible
hexagons and triangles of the Farey tesselation: a hexagon $W$ with vertices
$\pm(m,n)$, $\pm(p,q)$, and~$\pm(r,s)$ corresponds to an ideal triangle
$\Delta=(m/n,p/q,r/s)$. For example, the standard hexagon $W_0$ corresponds to
$\Delta_0=(\infty,-1,0)$. It is easy to check that a centrally symmetric
hexagon~$W$ is admissible if and only if $\Delta$ is one of the triangles that
form the Farey tesselation: both conditions are equivalent to the equations
$|mq-np|=|ps-qr|=|rn-sm|=1$. Moreover, flips of admissible hexagons correspond
to reflections in sides of the corresponding triangles. Consequently, the
graph~$\G$ introduced in~\S2.1 is nothing but the dual graph of the
tesselation; by the way, now it is obvious that~$\G$ is a tree.

To embed~$\G$ in~$H^2$ in the most symmetric way, let its vertices be the
centers of incircles\footnote{note that the center of a circle in~$H^2$ usually
does not coincide with its Euclidean center} of corresponding ideal triangles.
The incircles of two triangles having a common side touch this side at the same
point, because the whole picture is invariant under reflection in this side. By
a straightforward calculation, we obtain that the radius of the incircle of an
ideal triangle equals~$\frac{\ln3}2$. So the hyperbolic length of any edge of
$\G\subset H^2$ is equal to~$\ln3$. Angles between edges at any vertex are
equal to~$2\pi/3$ because of symmetry. This metric information completely
defines the embedding of~$\G$ in~$H^2$.

The embedding $\G\emb H^2$ described above (and represented on Fig.~22 by
dotted lines) makes the Farey tesselation coincide with the Voronoi partition
of~$H^2$ with respect to the vertices of~$\Gamma$: each ideal triangle~$\Delta(
w)$ containing a vertex $w\in\G$ coincides with the set of points $x\in H^2$
such that the vertex of~$\G$ closest to $x$ is~$w$. Moreover, for any point $x
\in\Delta(w)$ the point $y\in\G$ closest to~$x$ lies in~$\Delta(w)$, too. The
graph~$\G$ is the one-di\-men\-si\-onal skeleton of the Delaunay decomposition
of~$H^2$ dual to the Farey tesselation (considered as a Voronoi partition).

There are two different distance functions on~$\G$: the hyperbolic distance~$d_
H$ between its vertices and the ``counting'' distance $d_c$, which is equal to
the number of edges in the shortest path between vertices. The distance $d_c$
also can be defined as the number of the Farey tesselation lines that intersect
the geodesic line between two given vertices; here vertices can be replaced by
arbitrary points in the corresponding triangles. This implies the following
statement.

\proclaim{Theorem~16} For any point $z\in H^2$ not belonging to the Farey
tesselation lines, we have $d_c(\A z,z)=c(\A)$. 	\qed\endproclaim

\proclaim{Conjecture~4} There exists some natural construction of a metric
space $(\Sigma,d)$, where its elements $\sigma_t\in\Sigma$ are closely related
to triangles in the Farey tesselation, $d$ is closely related to~$d_c$, and
properties \rom{1--5} \rom(see the beginning of~\,\S\rom{4.1)} are satisfied.
\endproclaim

Unimodular quadratic forms also are natural candidates for the role of $\sigma_
t\in\Sigma$. For example, given a $\th$-curve, consider a quadratic form $Q_W$
on the lattice $\pi_1(T^2)$ (or rather on $H^1(T^2,\Z)$, which is the same)
that takes the vertices of the corresponding admissible hexagon~$W$ to~1. Then
$Q_W$ is a unimodular form with integer coefficients. Let us assign to $z\in H^
2$ an ellipse $E_z$ with semiaxes $e^t$ and $e^{-t}$, where $t=d_H(z,i)$, such
that the oriented angle between the $Ox$ axis and the large axis of $E_z$ is
half the oriented angle between the geodesic lines going from $i$ to zero and
to~$z$. Let $Q_z$ be a quadratic form defined by the condition $Q_z(\a)=1$ if
and only if $\a\in E_z$. Then $Q_W=Q_{z(W)}$, where $z(W)$ is the vertex
of~$\G$ that lies in the triangle $\Delta$ corresponding to an admissible
hexagon~$W$; this can be verified by a straightforward calculation.

\subhead 4.3. Other 3-dimensional manifolds \endsubhead

\subsubhead 4.3.1. Lens spaces \endsubsubhead We have already discussed lens
spaces in~\S1.3 and~\S2.4. Here we show that the distance function~$d_c$ coming
from the Farey tesselation leads to an elegant reformulation of Conjecture~1.
Note that $d_c$ is well-defined for pairs of rational points on the absolute
of~$H^2$ and for pairs where one point belongs to~$H^2$ and the other is a
rational absolute point.

\proclaim{Lemma~12} Let $p$, $q$ be coprime positive integers. Then the
Euclid complexity $E(p,q)$ defined in~\S\rom{1.2} is equal to $d_c(e^{2\pi i/3}
,p/q)$. 							\endproclaim

The point $e^{2\pi i/3}\in H^2$ can be replaced by any other inner point of the
triangle $\Delta_0=(\infty,-1,0)$.

\demo{Proof} By Theorem~5, we have $d(W,W_0)=E(p,q)$, where the admissible
hexagon $W$ is uniquely defined by its leading vertex~$(p,q)$ (see~\cite2 or
\S2.1 for definitions). Since $e^{2\pi i/3}\in\Delta_0$, we have $d_c(e^{2\pi i
/3},p/q)=d(W,W_0)$.						\qed\enddemo

Recall that any lens space $L_{p,q}$ is obtained by gluing together two solid
tori along their boundary torus~$T^2$. Denote by $\mu_1,\mu_2\in\pi_1(T^2)$
the meridians (that is, contractible boundary cycles) of two solid tori. Fix
some basis in $\pi_1(T^2)=\Z^2$. Then the slopes $r_1=r(\mu_1),r_2=r(\mu_2)\in
\Bbb Q$ of the cycles $\mu_1$ and~$\mu_2$ are defined.

\proclaim{Theorem~17} We have $d_c(r_1,r_2)=E(p,q)-1$.		\endproclaim

Compare this with Lemma~5 in~\S2.4.

\demo{Proof} The number $d_c(r(\mu_1),r(\mu_2))$ is independent of the choice
of a basis in~$\Z^2$, because the action~(4) of~$\SL(2,\Z)$ on~$H^2$ preserves
the Farey tesselation and, consequently, the function~$d_c$. For $L_{p,q}$,
there exists a basis of~$\Z^2$ such that $r_1=0$ and $r_2=p/q$. Since $p>q>0$,
we have $p/q\in(1,\infty)$. Hence, the distance $d_c$ of $p/q$ to~$0$ is one
less than its distance to an inner point of the triangle $(\infty,-1,0)$, see
Fig.~22. Now the statement of the theorem follows from Lemma~12.   \qed\enddemo

Thus Conjecture~1 is equivalent to the relation $c(L_{p,q})=d_c(r_1,r_2)-2$.

\remark{Remark} For any two (and even three) absolute points $a,b\in\R\cup
\infty$ there exists an isometry $A\in\SL(2,\R)$ of~$H^2$ taking them, say, to
0 and~1. This is no longer the case for the $\SL(2,\Z)$-action: obviously, only
rational points (including~$\infty$) can be mapped to rationals, but even a
pair of rational numbers $(r_1,r_2)$ cannot be mapped to the pair $(0,1)$
whenever $d_c(r_1,r_2)>0$. {\sl Exercises\/}: 1)~given two pairs of rational
numbers, determine whether they are $\SL(2,\Z)$-equi\-valent. (Hint:
see~\cite{23}; the condition $d_c(r_1,r_2)=d_c(s_1,s_2)$ is necessary but not
sufficient.) \,2)~Fix a basis in $\pi_1(T^2)$ and attach two solid tori along
$T^2$ so that their meridians have slopes $r_1$ and~$r_2$. We get a manifold $M
(r_1,r_2)$, which is a lens space. Show that $M(r_1,r_2)$ and $M(s_1,s_2)$ are
homeomorphic if and only if unordered pairs $(r_1,r_2)$ and $(s_1,s_2)$ are
$\SL(2,\Z)$-equi\-valent. 3)~Prove (once again) Theorem~3, see~\S1.3.
\endremark

\subsubhead 4.3.2. Stallings manifolds \endsubsubhead Consider a fibration $p\:
M^3\overset{F_g}\to\tto S^1$, where $F_g$ is an orientable surface (of
genus~$g$) and $M^3$ is an orientable 3-ma\-ni\-fold. It may happen that a
manifold~$M$ can be fibered over a circle in several different ways, and the
genus of the fiber need not be defined uniquely by~$M^3$. Choose any of the
fibering structures and denote by~$\A$ the monodromy, which is an isotopy class
of self-diffeo\-morphisms $A\:F_g\to F_g$.

The ideas of Section~3 can be applied to this situation as well. Nothing
changes in~\S3.1. Further, Lemma~7 remains true. Obviously, the set of isotopy
classes of \thcs should be replaced by the set of isotopy classes of trivalent
graphs $L\subset F_g$ such that $F_g\sm L$ is a 2-cell. With this correction,
analogues of Lemma~8 and Theorem~13 hold, while the proof of Lemma~8 requires
slight modification in the second paragraph. However, the difference between
1-skeletons $L_-$ and $L_+$ of $F_g$ in an analogue of Lemma~9 is measured by
at most one Dehn twist (rather than at most one flip), which may require up to
$4g-3$ flips. So in Theorem~14 we only get the estimate $d(W_+,W_-)\le5(4g-3)$,
where $W_+$, $W_-$ are isotopy classes of trivalent graphs embedded in~$F_g$ so
that their complements are 2-cells, and $d$ is the ``flip-distance''. The
remaining part of the reasoning of~\S3 goes smoothly, and we get an analogue of
Theorem~12 with $1/5(4g-3)$ instead of~$1/5$. The approach discussed in~\S4.1
can be applied to this situation, too.

Another difference from the case of $g=1$ is that we no longer know how to find
$c(\A)$, that is, the minimal number of flips required to convert a trivalent
graph $L$ such that $F_g\sm L$ is a 2-cell into its monodromy image~$\A L$.
Recall that $c(\A)$ can be computed by the methods of \S2.2 or~\S4.2 whenever
$g=1$. The former approach is obstructed by the fact that the graph~$\G$ is not
a tree if~$g>1$. Probably, some estimates for $c(\A)$ are easier to obtain than
its exact value. In order to do it using the latter approach, one should
replace the Farey tesselation of~$H^2$ by the cell decomposition of the
Teichm\"uller space of~$F_g$ described in~\cite4.

The construction of a spine with small number of vertices presented in~\S2.3
gives a spine of~$M^3$ with $c(\A)+4g+2$ vertices, and at least one of them may
be cancelled by simplification moves provided that $c(\A)>0$. It is not
immediately clear how many vertices can be cancelled out of $4g+2$. Anyway, it
is plausible that any spine of~$M^3$ contains more than $c(\A)$ vertices, where
$\A$ is the monodromy of any possible fibering of~$M^3$ over the circle.

\subsubhead 4.3.3. Topological economy principle \endsubsubhead If Conjecture~2
holds, a minimal spine of a torus bundle space $M(\A)$ can be constructed as
follows (see~\S2.3): fix a fiber $F=p^{-1}(0)\subset M^3$, then choose a \thc
$L\subset F$ that requires $c(\A)$ flips only to be converted to $\A L$,
construct a simple polyhedron~$P_0$ from the evolution of~$L$, and, finally,
add to~$P_0$ an extra face that cuts any path connecting two boundary
components of~$M^3\sm F$ in the complement of~$P_0$.

According to Thurston's geometrisation conjecture, any prime orientable compact
3-manifold $M^3$ can be cut into geometric pieces $M_i$, $i=1,\dots,m$, along
incompressible tori $T^2_j$, $j=1,\dots,n$, see, for example, the last section
of~\cite{22}. For nontrivial cycles in the $T^2_j$ remain nontrivial in~$M^3$
(by incompressibility), any spine $P$ of $M$ intersects all nontrivial cycles
in any~$T^2_j$; assuming general position, the intersections $P\cap T^2_j$
contain \thcs for all~$j$ by virtue of Lemma~8.

Given a family of \thcs $L_j\subset T^2_j$, consider simple polyhedra $P_i
\subset M_i$ such that $P_i\cap T^2_j=L_j$ for any boundary torus $T^2_j$ of~$M
_i$. Then add extra faces to~$P_i$ as necessary to get spines of the~$M_i$. The
union $P$ of obtained polyhedra is a spine of $M^3$ (maybe, with several
punctures). Minimize the total number of vertices of~$P$ over all possible
choices of the $P_i$ and the extra faces and over all ``boundary conditions'',
that is, all families $L_j\subset T^2_j$. Needless to say, it is not clear yet
how to implement this program in the general case.

\proclaim{Conjecture~5} A minimal spine of any \rom3-manifold \rom(that can be
cut into geometric pieces\rom) can be obtained by the procedure described in
the paragraph above.						\endproclaim

In other words, there exists a minimal spine of~$M^3$ that intersects any torus
$T^2_j$ along one \thc only. Similar conjecture can be formulated about
graph-mani\-folds~\cite{8,~30}, which can be cut along incompressible tori
into several copies of $D^2\x S^1$ and $N^2\x S^1$, where $N^2$ stands for $D^
2$ with two holes. Conjecture~5 is yet another facet of the ``topological
economy principle'' expressed and illustrated in~\cite3.

\enddocument
\end